\documentclass[10pt]{article}
\usepackage{graphicx}
\usepackage{epsfig}
\usepackage{amssymb}
\usepackage{amsmath}
\usepackage{amsfonts}
\usepackage{dsfont}
\newcommand{\RR}{\mathbb{R}}
\newcommand{\NN}{\mathbb{N}}
\newcommand{\ZZ}{\mathbb{Z}}

\def\div{\operatorname{div}}
\def\curl{\operatorname{curl}}

\newtheorem{theorem}{Theorem}[section]
\newtheorem{lemma}[theorem]{Lemma}
\newtheorem{corollary}[theorem]{Corollary}
\newtheorem{proposition}[theorem]{Proposition}
\newtheorem{e-definition}[theorem]{Definition}
\newtheorem{remark}{\bf Remark\/}

\def\longrightharpoonup{\relbar\joinrel\rightharpoonup}
\renewcommand{\div}{\operatorname{div}}

\begin{document}

\title{On correctors for linear elliptic homogenization in the presence of local defects}
\author{X. Blanc$^1$, C. Le Bris$^2$ \& P.-L. Lions $^3$\\
{\footnotesize $^1$ Universit\'e Paris Diderot, Laboratoire Jacques-Louis Lions,}\\
{\footnotesize B\^atiment Sophie Germain, 5, rue Thomas Mann}\\
{\footnotesize 75205 Paris Cedex 13, FRANCE,}\\
 {\footnotesize{\tt blanc@ann.jussieu.fr}}\\
{\footnotesize $^2$ Ecole des Ponts and INRIA,}\\
{\footnotesize 6 \& 8, avenue Blaise Pascal, 77455 Marne-La-Vall\'ee
  Cedex 2, FRANCE}\\
{\footnotesize{\tt lebris@cermics.enpc.fr}}\\
{\footnotesize $^3$ Coll\`ege de France, 11, place Marcelin Berthelot,}\\
{\footnotesize 75231 Paris Cedex 05,  and}\\
{\footnotesize CEREMADE, Universit\'e Paris Dauphine, Place de Lattre de
  Tassigny,}\\
{\footnotesize 75775 Paris Cedex 16, FRANCE}\\
{\footnotesize \tt lions@ceremade.dauphine.fr}
}
\maketitle

\begin{abstract} We consider the corrector equation associated, in homogenization theory, to a linear second-order
  elliptic equation in divergence form $-\partial_i(a_{ij}\partial_ju)=f$, when the diffusion coefficient is a
  locally  perturbed periodic coefficient. The question under study is the existence (and uniqueness) of the
  corrector, strictly sublinear  at infinity, with gradient in $L^r$ if the local perturbation is itself $L^r$,
  $r<+\infty$. The present work follows up on our works~\cite{milan,cras-defauts,cpde-defauts}, providing an
  alternative, more general and versatile approach, based on an {\it a priori} estimate, for this well-posedness
  result. Equations in non-divergence form such as~$-a_{ij}\partial_{ij}u=f$ are also considered, along with
  various extensions. The case of general advection-diffusion equations~$-a_{ij}\partial_{ij}u+b_j\partial_ju
  =f$ is postponed until our future work~\cite{BLL-2017-2}. An appendix contains a corrigendum to our
  earlier publication~\cite{cpde-defauts}.
 \end{abstract}
 \bigskip

\noindent{\bf Keywords:} Homogenization; elliptic PDE; $L^p$ estimates; defects.

\medskip

\noindent{\bf AMS subject classification:} 35J15, 35J70, 35B27, 74Q15, 76M50

\newpage
\tableofcontents
\newpage
\section{Introduction}

\paragraph{Motivation.} In a series of works~\cite{milan,cras-defauts,cpde-defauts} (see also the related works~\cite{josien,josien-these}), we have shown that the solution to a highly oscillatory equation of the type 
 \begin{equation}
  \label{eq:equation}
  -\hbox{\rm div}\, \left(a(x/\varepsilon)\,\nabla u^\varepsilon\right)=f
\end{equation}
may be efficiently approximated using the same ingredients as classical periodic homogenization theory when the coefficient $a$ in~\eqref{eq:equation}  is a \emph{perturbation} of a periodic coefficient, say to fix the ideas $a=a^{per} + \tilde a$ where $a^{per}$ is periodic and $\tilde a\in L^r$, $1\leq r<+\infty$, is a local perturbation that formally vanishes at infinity. The quality of the approximation (that is, the rate of convergence in $H^1$ norm of $u_\varepsilon$ minus its approximation based on homogenization theory) is entirely based upon the existence of a corrector function $w_p$, strictly sublinear at infinity (that is,  
$ \displaystyle \frac{ w_{p}(x)}{1+|x|}\quad
\mathop{\longrightarrow}^{|x|\to \infty}
\quad 0$), solution, for each $p\in \RR^d$, to the corrector equation associated to~\eqref{eq:equation}, namely
\begin{equation}
  \label{eq:correcteur}
  -\,\hbox{\rm div}\left(a\,(p+\nabla w_{p})\right)=0 \quad \text{in } \RR^d.
\end{equation}
Such a situation comes in sharp contrast to the general case of homogenization theory where only a sequence of "approximate" correctors~$w_{p,\varepsilon}$, satisfying $\displaystyle  -\,\hbox{\rm div}\left(a(x/\varepsilon)\,(p+\nabla w_{p,\varepsilon})\right)\quad\buildrel \varepsilon\rightarrow 0\over\rightharpoonup 0$, is needed to conclude, but where the rate of convergence of the approximation is then unknown. See more details in our previous works and in~\cite{josien,josien-these}. 
A quick inspection on~\eqref{eq:correcteur} shows that the corrector $w_p$ is expected to read as~$w_p=w_{p,per}+\tilde w_p$ with $w_{p,per}$ the periodic corrector (solution to $-\,\hbox{\rm div}\left(a^{per}\,(p+\nabla w_{p,per})\right)=0$) and $\tilde w_p$ solution with $\nabla \tilde w_p\in L^r$ to 
\begin{equation}
  \label{eq:correcteur-tilde}
  -\,\hbox{\rm div}\left(a\,\nabla \tilde w_{p}\right)= \,\hbox{\rm div}\left(\tilde a\,(p+\nabla w_{p,per})\right)\quad \text{in } \RR^d.
\end{equation}
In turn, the setting being linear, the existence and uniqueness of $\tilde w_p$ solution
to~\eqref{eq:correcteur-tilde} in the correct functional space is formally equivalent to the existence of an {\it a priori} estimate 
\begin{equation}
  \label{eq:estimee-informelle-div}
  \left\|\nabla u\right\|_{L^q}\,\leq C\, \left\|f\right\|_{L^q}, 
\end{equation}
for the exponent $q=r$, and $u$ solution to 
\begin{equation}
  \label{eq:equation-informelle-div}
  -\,\hbox{\rm div}\left(a\,\nabla u\right)= \,\hbox{\rm div}\,f\quad \text{in } \RR^d.
\end{equation}
The purpose of this article is to show how the estimate~\eqref{eq:estimee-informelle-div} (and similar estimates)
can be established with a good degree of generality (in particular $q$ needs not be equal to $r$, the Lebesgue
exponent such that $\tilde a\in L^r$, but can be any exponent $1\,< q\,<\,+\infty$), using a quite versatile
approach based on a simple version of the concentration-compactness principle~\cite{PLL-cc}. Intuitively,
estimate~\eqref{eq:estimee-informelle-div} holds true because the perturbation~$\tilde a$ within the
coefficient~$a$ in~\eqref{eq:equation-informelle-div} vanishes in a loose sense at infinity, while, by the
celebrated results of Avellaneda and Lin (see~\cite{AL1987,AL1989} and more specifically~\cite{AL1991}), the
estimate holds true when $a=a^{per}$. Thus, the integrability in $\RR^d$ of the solution remains unchanged. The
approach introduced here not only provides an alternate proof of the results of our earlier works for local
perturbations of periodic coefficients, but also allows for considering, instead of~\eqref{eq:equation}, equations
not in divergence form $\displaystyle  -\, a_{ij}\,\partial_{ij} u=f$, which were not approached in our works so
far. This also provides an approach for the case of advection-diffusion equations $\displaystyle  -\, a_{ij}(\cdot/\varepsilon)\,\partial_{ij} u^\varepsilon +\varepsilon^{-1}\,b_i(\cdot/\varepsilon)\,\partial_{i} u^\varepsilon =f$ which will be discussed in a forthcoming publication~\cite{BLL-2017-2}.

\medskip

\paragraph{Mathematical setting.} More precisely, ~\eqref{eq:equation} is supplied with homogeneous Dirichlet
boundary conditions and posed on a bounded regular domain $\Omega\subset\RR^d$,  with a right-hand-side term~$f\in
L^2(\Omega)$. For our exposition, we will assume $d\geq
2$. Of course, as always, dimension $1$ is specific and can be addressed by (mostly explicit) analytic arguments that we
omit here. 
In our earlier publications~\cite{milan,cras-defauts,cpde-defauts}, the coefficient~$a$ considered is  of the form~$a=a_{0}+\tilde a$ where $a_0$ denotes the unperturbed background, and $\tilde a$ the perturbation. For some of our results there, the unperturbed background can be quite general provided it enjoys the "natural" properties that make homogenization explicit. Similarly, we consider different cases of perturbations~$\tilde a$, and can prove some of the results in the absence of some regularity of the coefficients. We refer the reader
to~\cite{milan,cras-defauts,cpde-defauts} for all the precise settings and statements regarding the above informal claims. In the present contribution, however, we only address the case 
\begin{equation}
  \label{eq:aper+tildea}
  a=a^{per}+\tilde a
\end{equation}
where $a^{per}$ denotes a \emph{periodic} unperturbed background, and $\tilde a$ the perturbation, with 
\begin{equation}
\label{eq:hyp1}
\left\{
\begin{array}{l}
a^{per}(x)+\tilde a(x) \quad \hbox{\rm and}\quad a^{per}(x) \quad \hbox{\rm are both uniformly elliptic, in}\,x\in\RR^d, 
\\
a^{per}\in \left(L^\infty(\RR^d)\right)^{d\times d},
\\
  \tilde a \in \left(L^\infty(\RR^d)\cap L^r(\RR^d)\right)^{d\times d}, \quad \hbox{\rm for some }\quad 1\leq r<+\infty
\\
   a^{per},\, \tilde a\in \left(C^{0,\alpha}_{\rm unif}(\RR^d)\left(\RR^d\right)\right)^{d\times d} \quad\hbox{\rm for some}\quad \alpha>0,
\end{array}
\right.
\end{equation}

Note that, actually, in the above assumptions, the fact that $\tilde a$ is bounded is implied by the assumption
$\tilde a \in C^{0,\alpha}_{\rm unif}(\RR^d)\cap L^r$. We nevertheless state it as above to highlight the fact that $\tilde a \in L^q$ for
any $q\geq r$.  

The reason why we make the above set of assumptions~\eqref{eq:hyp1}  is that (a) we need our results to hold true
in the absence of the perturbation $\tilde a$  and the periodic case $a_0=a^{per}$ with $a^{per}$ sufficiently
regular is the only one where we are actually aware of (thanks to the works of Avellaneda and Lin) that this is the
case (see however Remark~\ref{rk:gloria} below), and (b) the perturbation~$\tilde a$ has to formally vanish at infinity for our specific arguments to hold.

We note, on the other hand, that we readily consider the case of matrix-valued coefficients, instead of scalar-valued
coefficients as in our previous works. The modifications are only a matter of technicalities. 

All the results of the present article are stated and proved for equation, not for systems. However, as we point
out in Remark~\ref{rk:systeme_2} below, the result of Proposition~\ref{prop:divergence-form} (i.e the divergence
form case) carries over to systems. This is not the case of our proof for the non-divergence form (see
Remark~\ref{rk:systeme_3} below), since our proof makes essential use of the maximum principle or of its consequences.

\medskip

 Given the above assumptions, it is well known \cite{blp} that there exists a periodic corrector $w_{p,per}$ unique up to the addition of a constant, that solves
 \begin{equation}
\label{eq:correcteur-per}
  -\,\hbox{\rm div}\left(a^{per}(x)\,(p+\nabla w_{p,per}(x))\right)=0,
\end{equation}
posed for each fixed vector $p\in \RR^d$. In these particular conditions,  $\nabla w_{p,per}\in C^{0,\alpha}_{\rm unif}(\RR^d)\cap L^\infty$. This corrector allows to consider the following first-order approximation to~$u^\varepsilon$ issued from the so-called two-scale expansion truncated at the first order
\begin{equation}
\label{eq:ueps1-per}
  u^{\varepsilon,1}_{per}(x)=u^*(x)+\varepsilon \sum_{i=1}^d\partial_{x_i}u^*(x)\,w_{e_i,per}(x/\varepsilon),
\end{equation}
where $e_i$ are the canonical vectors of $\RR^d$ and where $u^*$ denotes the homogenized limit of~$u^\varepsilon$, that is the solution to 
\begin{equation}
\label{eq:equation*}
-\hbox{\rm div}\, \left(a^*\,\nabla u^*\right)=f,
\end{equation}
with homogeneous
Dirichlet boundary conditions on $\partial\Omega$, where $a^*$ is the homogenized matrix-valued coefficient (actually also computed from local averages of the solution~$w_{p,per}$ to~(\ref{eq:correcteur-per})). We have that $u^\varepsilon -u^{\varepsilon,1}_{per}$ converges to zero (at least) in $H^1$ norm  and, precisely because of the existence of the corrector, the rate of the convergence $\displaystyle \left \|u^\varepsilon -u^{\varepsilon,1}_{per}\right \|_{H^1}$ as $\varepsilon$ vanishes
can be made precise in terms of $\varepsilon$. We refer the reader to our previous works and the classical textbooks~\cite{blp,jikov} for more details. 

It has been pointed out in our works that this quality of approximation carries over to the case of a local
perturbation of the coefficient in~\eqref{eq:equation}. A proof of this fact is contained in
\cite{josien,josien-these}. Problem~(\ref{eq:correcteur-per}), now reading
as~\eqref{eq:correcteur},   is therefore a key step in the understanding, and approximation, of the
solution~$u^\varepsilon$ both locally and globally. This fact is intuitively clear when one has realized that this
problem is obtained by zooming-in from~(\ref{eq:equation}) to the small scale. Using linearity, ~\eqref{eq:correcteur} is equivalent to~\eqref{eq:correcteur-tilde} and the key question is thus to prove existence for the latter equation.

\medskip

\paragraph{Plan.} Our contribution is organized as follows. To start with, we consider in Section~\ref{sec:divergence} the case of the equation in divergence form~\eqref{eq:equation} under the conditions made precise in~\eqref{eq:hyp1}. We establish the estimate announced in~\eqref{eq:estimee-informelle-div} for the solutions to~\eqref{eq:equation-informelle-div}. Our result is stated in Proposition~\ref{prop:divergence-form}. The subsequent Section~\ref{sec:non-divergence} is devoted to the analogous estimate, stated in Proposition~\ref{prop:non-divergence-form}, for the equation in non-divergence form. The fact that each of the estimate implies the well-posedness of the corresponding corrector problem (and thus, subsequently and using the arguments of our other works, the agreement
 of the first order approximation~\eqref{eq:ueps1-per} with the oscillatory solution~$u^\varepsilon$ in various
 norms and at a certain well defined rate in~$\varepsilon$) is made precise in
 Section~\ref{sec:homogenization}. Finally, we take the opportunity of the present article to provide, in
 Appendix~\ref{sec:Appendix}, a corrigendum of our previous work~\cite{cras-defauts,cpde-defauts}. Although
 this did not at all affect the main results of our works, we made there, for some  intermediate technical result
 (namely Lemma 4.2 of~\cite{cpde-defauts} and Lemma~1 of \cite{cras-defauts}), some erroneous claims.  We correct this here.

\section{Estimate for operators in divergence form}
\label{sec:divergence}
%
%
As mentioned in our introduction, we wish to prove existence and uniqueness of the (strictly sublinear at infinity)
corrector $w_p$, solution for $p\in\RR^d$ fixed, to the corrector equation~\eqref{eq:correcteur}. Assuming that the
coefficient~$a$ is of the form~\eqref{eq:aper+tildea} and satisfies the assumptions~\eqref{eq:hyp1}, we readily
introduce $\tilde w_p=w_p-w_{p,per}$ where the latter denotes the periodic corrector associated to $a^{per}$, the
existence and uniqueness (up to the addition of a constant) of which is a classical fact \cite{blp}. For further
reference, we note that under the regularity conditions~\eqref{eq:hyp1} for the coefficient~$a^{per}$, elliptic
regularity implies that the periodic corrector~$w_{p,per}$ is a $W^{1,\infty}$ function. Indeed, the
classical Hilbert theory gives that $w_{p,per}\in H^1_{loc}$ and is periodic. Harnack inequality then implies that
$w_{p,per}\in L^\infty$. Finally, \cite[Theorem~8.32]{GT} implies that $\nabla w_{p,per}$ is H\"older continuous,
hence in particular is in $W^{1,\infty}$. Equation~\eqref{eq:correcteur} reads as~\eqref{eq:correcteur-tilde}, that is,
$$-\,\hbox{\rm div}\left(a\,\nabla \tilde w_{p}\right)= \,\hbox{\rm div}\left(\tilde a\,(p+\nabla w_{p,per})\right)\quad \text{in } \RR^d,
$$
which we now have to solve. Formally simplifying both sides of the equation leads to considering the
equation~$-\Delta \tilde w_{p}=\,\hbox{\rm div}\left(\tilde a\,p\right)$ and we thus expect, for $r>1$, to find $\nabla \tilde w_p$ in the same space as $\tilde a$, namely $L^r$. This of course will in particular ensure that $w_p$ is strictly sublinear at infinity. This expectation is confirmed by the results of our earlier contributions, which we now prove in a different manner here. Our main result is the following:

\begin{proposition}
\label{prop:divergence-form}
Assume \eqref{eq:aper+tildea}-\eqref{eq:hyp1}. Fix $1<q<+\infty$. Then, for all $f\in \left(L^q(\RR^d)\right)^d$, there exists~$u\in L^1_{loc}(\RR^d)$, such that $\nabla u\in \left(L^q(\RR^d)\right)^d$, solution to equation~\eqref{eq:equation-informelle-div} namely
\begin{equation}
  \tag{\ref{eq:equation-informelle-div}}
-\,\hbox{\rm div}\left(a\,\nabla u\right)= \,\hbox{\rm div}\,f\quad \text{in } \RR^d.
\end{equation}
Such a solution is  unique up to the addition of a constant. In addition, there exists a constant $C_q$, independent on $f$ and $u$, and only depending  on $q$, $d$ and the coefficient $a$, such that   $u$ satisfies~\eqref{eq:estimee-informelle-div}, namely
\begin{equation}
    \left\|\nabla u\right\|_{\left(L^q(\RR^d)\right)^{d}}\,\leq C_q\, \left\|f\right\|_{\left(L^q(\RR^d)\right)^{d}}. \tag{\ref{eq:estimee-informelle-div}}
\end{equation}
\end{proposition}

\medskip

The existence and uniqueness (up to the addition of a constant) of~$\tilde w_p$ (and thus of the corrector~$w_p$)
is an immediate consequence of Proposition~\ref{prop:divergence-form}. For $r>1$, the proposition is applied to $q=r$, $f=\tilde a\,(p+\nabla w_{p,per})$, given
that $\tilde a\in (L^r(\RR^d))^{d\times d}$ and~$\nabla w_{p,per}\in (L^\infty(\RR^d))^d$. The case $r=1$ is
considered in Remark~\ref{rk:notL1} below.

\begin{remark}\label{rk:notL1}
As stated in assumption \eqref{eq:hyp1}, the case $\tilde a \in L^r$ for
$r=1$ is allowed in Proposition~\ref{prop:divergence-form}. However, the proposition then only gives existence of $\nabla
\tilde
w_p\in L^q(\RR^d)$ for any $q>1$, and not $\nabla \tilde w_p\in L^1$ as $\tilde a$. Writing $-\div(a^{per}\nabla \tilde w_p) = \div(\tilde a (\nabla w_p+p))$ and
using $\nabla w_p\in L^\infty$, a fact that is established there, the
results of \cite{AL1991} imply that $\nabla \tilde w_p$ is in weak-$L^1(\RR^d)$. But $\nabla \tilde w_p\notin
L^1(\RR^d)$, contrary to what is mistakenly stated in \cite{cpde-defauts,cras-defauts}. A counterexample for $a^{per} = 1,$ $\tilde a$
with compact support, is provided in
Remark~\ref{rk:L1} below.
\end{remark}

\medskip

The proof of Proposition~\ref{prop:divergence-form} to which we now proceed is based upon the following intuitive
property. When the defect~$\tilde a$ identically vanishes, the coefficient~$a$ is the periodic coefficient~$a^{per}$. In this particular case, estimate~\eqref{eq:estimee-informelle-div} has been established in~\cite{AL1991}. The estimate is shown there, using the representation of the solution $u$ in terms of  the Green function $G^{per}(x,y)$ associated to the operator~$-\,\hbox{\rm div}\left(a^{per}\,\nabla .\right)$, and the properties of approximation of this Green function obtained from the results of~\cite{AL1987}. Next, when~$\tilde a\not=0$, one notices that, since~$\tilde a\in \left( L^r(\RR^d)\cap C^{0,\alpha}_{\rm unif}(\RR^d)(\RR^d)\right)^{d\times d}$, we have that $\tilde a(x)\buildrel |x|\rightarrow \infty\over\longrightarrow 0$. Intuitively, the operator $-\,\hbox{\rm div}\left(a\,\nabla .\right)$ is therefore close to the operator~$-\,\hbox{\rm div}\left(a^{per}\,\nabla .\right)$  at infinity, and estimate~\eqref{eq:estimee-informelle-div} is likely to hold true there. On the other hand, locally, estimate~\eqref{eq:estimee-informelle-div} is a consequence of elliptic regularity and the fact that it holds true in the Hilbertian case $q=2$. The actual rigorous proof of Proposition~\ref{prop:divergence-form} implements this strategy of proof, using a continuation argument, the celebrated results of~\cite{AL1991} and our results~\cite{milan} on the case $q=2$.

\medskip

\noindent{\bf Proof of Proposition~\ref{prop:divergence-form}} We argue by continuation. We henceforth fix some~$2\leq q<+\infty$. The case $1<q\leq 2$ will be obtained by duality at the end of the proof. 

We define~$a_t=a^{per}+t\,\tilde a$ and intend to prove the statements of Proposition~\ref{prop:divergence-form} for $t=1$. For this purpose, we introduce the property ${\mathcal P}$ defined by: {\it we say that the coefficient $a$, satisfying the assumptions~\eqref{eq:aper+tildea}-\eqref{eq:hyp1} (for some $1\leq r<+\infty$) satisfies ${\mathcal P}$ if the statements of Proposition~\ref{prop:divergence-form} hold true for equation~\eqref{eq:equation-informelle-div} with coefficient~$a$}. We next define the interval
\begin{equation}
\label{eq:interval}
{\mathcal I}= \left\{t\in [0,1]\,/\, \forall\, s\in [0,t], \hbox{\rm Property}\,{\mathcal P}\,\hbox{\rm is true for}\, a_s\right\}.
\end{equation}
We intend to successively prove that ${\mathcal I}$ is not empty, open and closed (both notions being understood relatively to the closed interval $[0,1]$), which will show that ${\mathcal I}=[0,1]$, and thus the result claimed.

\medskip

To start with, we remark that the results of Avellaneda and Lin in~\cite[Theorem A]{AL1991} show that ${\mathcal I}\not=\emptyset$ since $0\in {\mathcal I}$. Notice that the property $u\in L^1_{loc}$ is a straightforward consequence of elliptic regularity using $f\in \left(L^1_{loc}(\RR^d)\right)^d$ and the H\"older regularity of the coefficient~$a^{per}$ because of~\eqref{eq:hyp1}. This property immediately carries over to  all~$t\in [0,1]$ as soon we know there is a solution in the following argument.

\medskip

Next, we show that ${\mathcal I}$ is open (relatively to the interval~$[0,1]$). For this purpose, we suppose that
$t\in {\mathcal I}$ and wish to prove property ${\mathcal P}$ on $[t,t+\varepsilon[$ for some $\varepsilon>0$. In order to solve, for $f\in \left(L^q(\RR^d)\right)^d$ 
$$ -\,\hbox{\rm div}\left((a_t+\varepsilon\,\tilde a)\,\nabla u\right)= \,\hbox{\rm div}\,f\quad \text{in } \RR^d,$$
we write it as follows:
\begin{equation}\label{eq:1}
  \nabla u = \Phi_t\left( \varepsilon \tilde a \nabla u + f\right),
\end{equation}
where $\Phi_t$ is the linear map which to $f\in L^q(\RR^d)$ associates $\nabla u\in L^q(\RR^d)$, where $u$ is the solution to
\eqref{eq:equation-informelle-div}. Since $\Phi_t$ is continuous from $L^q$ to $L^q$, with norm $C_q$, it is clear
that, for $C_q \varepsilon \left\|\tilde a\right\|_{L^\infty(\RR^d)} <1$, the above map is a contraction. Hence,
applying the Banach fixed-point theorem, \eqref{eq:1} has a unique solution in $L^q(\RR^d)$, which satisfies the estimate 
\eqref{eq:estimee-informelle-div}, in which $C_q$ has been replaced by $C_q\left(C_q \varepsilon \left\|\tilde a\right\|_{L^\infty(\RR^d)}-1\right)^{-1}.$
\medskip

We now show, and this is the key point of the proof,  that ${\mathcal I}$ is closed. We assume that $t_n\in
{\mathcal I}$, $t_n\leq t,$ $t_n\longrightarrow t$ as $n\longrightarrow +\infty$. For all $n\in\NN$ we know that, for any $f\in \left(L^{q}(\RR^d)\right)^{d}$, we have a solution (unique to the addition of a constant) $u$ with~$\nabla u\in \left(L^{q}(\RR^d)\right)^{d}$ of the equation 
$$-\,\hbox{\rm div}\left(a_{t_n}\,\nabla u\right)= \,\hbox{\rm div}\,f\quad \text{in } \RR^d,
$$
and that this solution satisfies $\left\|\nabla u\right\|_{\left(L^q(\RR^d)\right)^{d}}\leq C_n\,\left\|f\right\|_{\left(L^q(\RR^d)\right)^{d}}$ for a constant $C_n$ depending on $n$ but not on~$f$ nor on~$u$. We want to show the same properties for $t$. 

We first temporarily admit that the sequence of constants $C_n$ is uniformly bounded from above in~$n$ and conclude.
For $f\in \left(L^{q}(\RR^d)\right)^{d}$ fixed, we consider the sequence of solutions~$u^n$ to $$-\,\hbox{\rm div}\left(a_{t_n}\,\nabla u^n\right)= \,\hbox{\rm div}\,f\quad \text{in } \RR^d,
$$
which we may write as
$$-\,\hbox{\rm div}\left(a_{t}\,\nabla u^n\right)= \,\hbox{\rm div}\,\left(f+(a_{t_n}-a_t)\,\nabla u^n\right)\quad \text{in } \RR^d,
$$

The sequence of gradients $\nabla u^{n}$ is bounded in~$\left(L^q(\RR^d)\right)^{d}$, and therefore weakly
converges (up to an extraction) to some $\nabla u$. We may pass  to the weak limit in the above equation (recall
that $a_{t_n}-a_t$ converges strongly in~$L^\infty$) and we find a solution to $-\,\hbox{\rm
  div}\left(a_{t}\,\nabla u\right)= \,\hbox{\rm div}\,f$. The solution also satisfies the estimate (because the
sequence $C_n$ is bounded and because the norm is weakly lower semi continuous). There remains to prove uniqueness,
that is, 
\begin{equation}
\label{eq:zero-1}
-\,\hbox{\rm div}\left(a\,\nabla u\right)= \,0\quad \text{in } \RR^d,
\end{equation}
with $\nabla u \in L^q(\RR^d)$,
implies $\nabla u\equiv 0$ in the present setting.
To this end, we notice that~\eqref{eq:zero-1} also reads as 
\begin{equation}
\label{eq:zero-1-bis}
-\,\hbox{\rm div}\left(a^{per}\,\nabla u\right)= \,t\,\hbox{\rm div}\left(\tilde a\,\nabla u\right).
\end{equation}
Using that, in the right-hand side,  $\tilde a\,\nabla u\in \left(L^{q_1}(\RR^d)\right)^{d}$ for $\displaystyle {1\over q_1}= {1\over r}+{1\over q}$ (by the H\"older inequality), and the results of~\cite{AL1991} on the operator with periodic coefficient, this implies that, in turn, $\nabla u\in \left(L^{q_1}(\RR^d)\right)^{d}$. 
One may then iterate this argument, and inductively obtain $\nabla u\in \left(L^{q_n}(\RR^d)\right)^{d}$ for~$\displaystyle {1\over {q_n}}= {1\over {q_{n-1}}}+{1\over r}$. 
One thereby obtains~$\nabla u\in \left(L^{q_n}(\RR^d)\right)^{d}$ for~$\displaystyle {1\over {q_n}}= {1\over
  q}+{n\over r}$. We recall that we have assumed $q\geq 2$, thus $\frac 1 q \leq \frac 1 2$. If in addition $r\geq
2$, it is then always possible to find $n\geq 0$ such that $1\leq q_n \leq 2$. 
In the case $r<2$,
we note that $\tilde a \in L^r\cap L^\infty \subset L^2$, hence we apply the argument with $r=2$. 
In both cases, we have $1\leq q_n\leq 2$ for some adequate $n\geq 0$, and we obtain, by interpolation, $\nabla u\in \left(L^{q_n}\cap L^{q}(\RR^d)\right)^{d}\subset
\left(L^{2}(\RR^d)\right)^{d}$. But, for such an $L^2$ function, \eqref{eq:zero-1} immediately implies $\nabla
u\equiv 0$ by ellipticity (a precise argument may be found in \cite{milan}). This concludes the argument of uniqueness.

\medskip

In order to prove that the sequence of constants $C_n$ is indeed bounded, we argue by contradiction and assume that the sequence $C_n$ is unbounded, which amounts to assuming there exist~$f^n\in \left(L^{q}(\RR^d)\right)^{d}$ and $u^n$ with~$\nabla u_n\in \left(L^{q}(\RR^d)\right)^{d}$, such that 
\begin{equation}
\label{eq:contradic1}
-\,\hbox{\rm div}\left(a_{t_n}\,\nabla u^n\right)= \,\hbox{\rm div}\,f^n\quad \text{in } \RR^d,\end{equation}
\begin{equation}
\label{eq:contradic2}
\left\|f^n\right\|_{\left(L^q(\RR^d)\right)^{d}}\buildrel n\longrightarrow +\infty\over \longrightarrow 0,
\end{equation}
\begin{equation}
\label{eq:contradic3}
\left\|\nabla u^{n}\right\|_{\left(L^q(\RR^d)\right)^{d}}=1,\quad\hbox{\rm for all}\,n\in \NN.
\end{equation}
We readily notice that~\eqref{eq:contradic1} also reads as 
$$-\,\hbox{\rm div}\left(a_{t}\,\nabla u^n\right)= \,\hbox{\rm div}\,\left(f^n+(a_{t}-a_{t_n})\,\nabla u^n\right),$$
where,  as $n\longrightarrow 0$, the rightmost term inside the divergence strongly converges to zero in $\left(L^q(\RR^d)\right)^{d}$ because $a_{t}-a_{t_n}$ strongly converges to zero in $\left(L^q(\RR^d)\right)^{d}$ and $\nabla u^n$ is bounded in $\left(L^q(\RR^d)\right)^{d}$. Therefore, without loss of generality, we may change $f^n$  into~$f^n+(a_{t}-a_{t_n})\,\nabla u^n$ and  replace~\eqref{eq:contradic1} by
\begin{equation}
\label{eq:contradic1-bis}
-\,\hbox{\rm div}\left(a_{t}\,\nabla u^n\right)= \,\hbox{\rm div}\,f^n\quad \text{in } \RR^d.\end{equation}
We now concentrate our attention on the sequence~$\nabla u^n$. In the spirit of the method of concentration-compactness, we now claim that 
\begin{equation}
\label{eq:cc-1}
\exists\, \eta>0, \quad \exists\, 0<R<+\infty,\quad \forall \,n\in\NN,\quad \left\|\nabla u^{n}\right\|_{\left(L^q(B_R)\right)^{d}}\geq \eta>0,
\end{equation}
where $B_R$ of course denotes the ball of radius~$R$ centered at the origin. 

We again argue by contradiction (we recall the main argument we are conducting here is also an argument by contradiction) and assume that, contrary to~\eqref{eq:cc-1},
\begin{equation}
\label{eq:cc-1-not}
\forall\, 0<R<+\infty,\quad  \left\|\nabla u^{n}\right\|_{\left(L^q(B_R)\right)^{d}}\buildrel n\longrightarrow +\infty\over \longrightarrow 0.
\end{equation}
Since $\tilde a$ satisfies the properties in~\eqref{eq:hyp1}, it vanishes at infinity and thus, for any~$\delta>0$,  we may find some sufficiently large radius~$R$ such that 
\begin{equation}
\label{eq:cc-2}
\left\|\tilde a\right\|_{\left(L^\infty(B_R^c)\right)^{d}}\leq \delta,
\end{equation}
where~$B_R^c$ denotes the complement set of the ball $B_R$.  We then write
\begin{eqnarray}
\left\|\tilde a\,\nabla u^n\right\|^q_{\left(L^q(\RR^d)\right)^{d}}&=&\int_{B_R}\left|\tilde a\,\nabla u^n\right|^q+\int_{B_R^c}\left|\tilde a\,\nabla u^n\right|^q\nonumber\\
&\leq&\left\|\tilde a\right\|^q_{\left(L^\infty(\RR^d)\right)^{d\times d}}\,\left\|\nabla u^n\right\|^q_{\left(L^q(B_R)\right)^{d}}\nonumber\\
&&+\left\|\tilde a\right\|^q_{\left(L^\infty(B_R^c)\right)^{d\times d}}\,\left\|\nabla u^n\right\|^q_{\left(L^q(\RR^d)\right)^{d}}\nonumber\\
&\leq&\left\|\tilde a\right\|^q_{\left(L^\infty(\RR^d)\right)^{d\times d}}\,\left\|\nabla u^n\right\|^q_{\left(L^q(B_R)\right)^{d}}\nonumber\\
&&+\delta\,.1,
\end{eqnarray}
using~\eqref{eq:contradic3} and~\eqref{eq:cc-2} for the latter majoration. 
On the other hand, ~\eqref{eq:cc-1-not} implies that the first term in the right hand side vanishes, and since
$\delta$ is arbitrary, this shows that $\tilde a\,\nabla u^n$ converges strongly to zero in~$\left(L^q(\RR^d)\right)^{d}$. Inserting this and~\eqref{eq:contradic2} into~\eqref{eq:contradic1-bis}, which, for this specific purpose, we rewrite as 
$$-\,\hbox{\rm div}\left(a^{per}\,\nabla u^n\right)= \,\hbox{\rm div}\,\left(f^n+t\,\tilde a\,\nabla u^n\right),
$$
and using the continuity result in~$\left(L^q(\RR^d)\right)^{d}$ for the periodic setting established by Avellaneda and Lin in~\cite{AL1991}, we deduce that $\nabla u^n$ (strongly) converges to zero in~$\left(L^q(\RR^d)\right)^{d}$. This evidently contradicts~\eqref{eq:contradic3} . We therefore have established~\eqref{eq:cc-1}.

Because of the bound~\eqref{eq:contradic3}, we may claim that, up to an extraction, $\nabla u^n$ weakly converges
in $\left(L^q(\RR^d)\right)^{d}$, to some $\nabla u$. Passing to the limit in the equation in the sense of distributions, we have
$-\div(a_t\nabla u) = 0$. Next, we show that this convergence in $\left(L^q_{loc}(\RR^d)\right)^{d}$ is indeed
strong. By Sobolev compact embeddings, we know this convergence implies the
strong convergence of the sequence $u^n$ to~$u$ (up to a sequence of irrelevant constants $c_n$ which, with a
slight abuse of notation, we may include in $u$)
in~$L^q_{loc}(\RR^d)$. Since we have
\begin{displaymath}
  -\div\left( a_t\left(\nabla u^n - \nabla u\right)\right) = \div(f^n),
\end{displaymath}
we multiply this equation by $(u^n-u)\chi_R$, where $\chi_R$ is a smooth cut-off function such that $\chi_R = 1$ in $B_R$,
and $\chi_R = 0$ in $B_{R+1}^c$. Integrating by parts, we have
\begin{multline*}
  \int  \chi_R\left[a_t\nabla (u^n-u)\right]\cdot \nabla(u^n-u) = -\int \left[a_t\nabla
      (u^n-u)\right]\cdot \left(\nabla \chi_R\right) (u^n-u) \\ + \int  \chi_R\left[\nabla (u^n-u)\right]\cdot
  f^n  + \int (u^n-u)\, f^n\cdot \nabla \chi_R .
\end{multline*}
We may pass to the limit in each term of the right hand side, since $\nabla (u^n - u) \longrightharpoonup 0$ in
$L^q$, hence in $L^2(B_R)$, $u^n-u$ converges strongly in $L^2(B_R)$, and $f^n\longrightarrow 0$ in
$L^q$, hence in $L^2(B_R)$. Using the fact that $a$ is elliptic, we therefore obtain that $\nabla u^n\longrightarrow \nabla u$
in $L^2(B_R)$. Finally, using the elliptic estimate (see e.g.~\cite[Theorem 7.2]{Giaquinta})
\begin{equation}\label{eq:20}
\left\|\nabla v\right\|_{\left(L^q(B_R)\right)^{d}}\leq\, C(R)\,\left(\left\|\nabla v\right\|_{\left(L^2(B_{2R})\right)^{d}}\,+\,\left\|f\right\|_{\left(L^q(B_{2R})\right)^{d}}\right),
\end{equation}
for all $0<R<+\infty$, and all solutions~$v$ to $-\,\hbox{\rm div}\left(a_t\,\nabla v\right)=\,\hbox{\rm div} f$,
we obtain, applying~\eqref{eq:20} to $u^n-u$, that $\nabla u^n$ strongly converges in $\left(L^q(B_R)\right)^{d}$
to  $\nabla u$. We infer from~\eqref{eq:cc-1} and the local strong convergence that $\nabla u$ cannot identically
vanish, while it solves $-\div(a_t\nabla u) = 0$. The argument following \eqref{eq:zero-1} implies that
$\nabla u = 0$, and we reach a final contradiction. This concludes the proof in the case $2\leq q<+\infty$. 

\medskip

For the case $1<q< 2$, we argue as announced by duality. At this stage, we have established the claims within Proposition~\ref{prop:divergence-form} for $2\leq q<+\infty$.  
We may apply them to the case of the transposed coefficient $a^T$ of $a$, and the operator $-\,\hbox{\rm div}\left(a^T\,\nabla .\right)$, since obviously, $a^T$ satisfies the assumptions~\eqref{eq:hyp1} if $a$ does. 
Let us now fix $f\in\left(L^q(\RR^d)\right)^{d}$, for some~$1<q\leq 2$, and denote by $2\leq q'<+\infty$, the conjugate exponent of $q$, that is, $\displaystyle {1\over q}+{1\over {q'}}=1$. To any arbitrary function~$g\in\left(L^{q'}(\RR^d)\right)^{d}$, we may associate the unique (up to an additive constant) solution~$v$, such that $\nabla v\in\left(L^{q'}(\RR^d)\right)^{d}$, 
to $-\,\hbox{\rm div}\left(a^T\,\nabla v\right)=\,\hbox{\rm div}\, g$. Its gradient $\nabla v$ depends linearly, continuously, on $g$. We may therefore define by $g\longmapsto L_f(g):=\int_{\RR^d} f\,.\,\nabla v$ a linear form
on $\left(L^{q'}(\RR^d)\right)^{d}$. Since, 
using the result of Proposition~\ref{prop:divergence-form} for $q'$,
\begin{eqnarray}
\label{eq:linear-form-1}
\left | L_f(g)=\int_{\RR^d} f\,.\,\nabla v\right |& \leq &\,\left\| f\right\|_{\left(L^q(\RR^d)\right)^d}\,
\left\| \nabla v\right\|_{\left(L^{q'}(\RR^d)\right)^d}\nonumber\\
&\leq& C_{q'}\,\left\| f\right\|_{\left(L^q(\RR^d)\right)^d}\,
\left\| g\right\|_{\left(L^{q'}(\RR^d)\right)^d},
\end{eqnarray}
 this linear form is therefore a continuous map on $\left(L^{q'}(\RR^d)\right)^d$. Hence there exists some $U\in \left(L^q(\RR^d)\right)^d$ such that 
 $$L_f(g)=\int_{\RR^d} f\,.\,\nabla v=\int_{\RR^d} g\,.\,U,$$
and we read on estimate~\eqref{eq:linear-form-1} that 
$$\left\|U\right\|_{\left(L^q(\RR^d)\right)^{d}}\leq C_{q'}\,\left\| f\right\|_{\left(L^q(\RR^d)\right)^d}.
$$
We now identify more precisely~$U$. Assuming that $g$ additionally satisfies $\hbox{\rm curl}\, g=0$ on $\RR^d$, we have $-\,\hbox{\rm div}\left(a^T\,\nabla v\right)=\,\hbox{\rm div}\, g=0$ with $\nabla v\in \left(L^{q'}(\RR^d)\right)^d$, and thus, by the estimate, $\nabla v\equiv 0$. It follows that, for such $g$, $L_f(g)=\int_{\RR^d} f\,.\,\nabla v=0$, thus $\int_{\RR^d} g\,.\,U=0$. This property shows that there exists some $u$ such that $U=\nabla u$, and thus $\nabla u\in \left(L^{q}(\RR^d)\right)^d$ with
$$\left\|\nabla u\right\|_{\left(L^q(\RR^d)\right)^{d}}\leq C_{q'}\,\left\| f\right\|_{\left(L^q(\RR^d)\right)^d}.
$$

We finally show that $u$ satisfies $-\,\hbox{\rm div}\left(a\,\nabla u\right)=\,\hbox{\rm div}\, f$. To this end,
we consider the specific case where $v\in \mathcal{D}(\RR^d)$  (that is, $v$ is smooth and has compact support) and set $g=a^T\,\nabla v$, so that, in effect, $-\,\hbox{\rm div}\left(a^T\,\nabla v\right)=\,\hbox{\rm div}\, g$ holds true. Applying the above, we have  $\displaystyle \int_{\RR^d} f\,.\,\nabla v=\int_{\RR^d} g\,.\,\nabla u$. The left-hand side is the duality  product $\langle\,\hbox{\rm div}\, f,v\rangle_{{\mathcal D'}(\RR^d),\mathcal{D}(\RR^d)}$, while, by definition of~$g$, the right-hand side reads as 
$$ \int_{\RR^d} a^T\,\nabla v\,.\,\nabla u=\int_{\RR^d} \,\nabla v\,.\,a\,\nabla u=\langle\,-\hbox{\rm div}\,\left(a\,\nabla u\right),v\rangle_{{\mathcal D'}(\RR^d),\mathcal{D}(\RR^d)}.$$
 Since this holds true for all~$v\in \mathcal{D}(\RR^d)$, this shows  $-\,\hbox{\rm div}\left(a\,\nabla u\right)=\,\hbox{\rm div}\, f$ and concludes our proof. \hfill$\diamondsuit$

\medskip
\begin{remark}\label{rk:systeme_1}
  Although, in the above proof, the case $1<q<2$ is proved by duality, it is also possible to use a direct approach
  similar to the above argument. The heart of the above proof is the following result: if $u$ is a solution to
  $-\div(a\nabla u) = 0$ and if $\nabla u \in L^q(\RR^d)$, then $u$ is constant. We use it for $q=2$, but it is
  still valid for any $q>1$, as it is stated in
  Lemma~\ref{lm:q<2} below. Note however that we do not know if this lemma, which is proved only for
  equations, carries over to systems. 
\end{remark}
\begin{lemma}\label{lm:q<2}
  Assume $d\geq 2$, and that the matrix $a$ is $C^{0,\alpha}_{\rm unif}(\RR^d)$, uniformly elliptic and bounded. If $u$ satisfies $-\div(a\nabla u) = 0$
  and $\nabla u\in L^q(\RR^d)$, for some $1\leq q <d$, then $\nabla u =0$. 
\end{lemma}
\noindent{\bf Proof:} We first note that, according to the Galiardo-Nirenberg-Sobolev inequality \cite[Section 5.6.1,
Theorem 1]{evans}, $u\in L^{q^*}(\RR^d)$, up to the addition of a constant, with $\frac 1 {q^*} = \frac 1 q - \frac 1 d$. Hence, we have $u\in
W^{1,q^*}_{\rm loc}$, with
\begin{displaymath}
  \sup_{x_0\in\RR^d} \|u\|_{W^{1,q^*}(B_2(x_0))} <+\infty.
\end{displaymath}
Here, $B_2(x_0)$ is the ball of radius $2$ centered at $x_0$. Next, we
apply \cite[Theorem~1]{brezis-2008}, which states that, if $u$ satisfies the above
properties, then for any ball $B_1(x_0)$ of radius $1$, we have $u\in W^{1,s}(B_1(x_0))$ for all $1<s<+\infty$, with the following estimate
\begin{equation}
  \label{eq:estimation_brezis}
  \|u\|_{W^{1,s}(B_1(x_0))} \leq C \|u\|_{W^{1,q^*}(B_2(x_0))},
\end{equation}
where $C$ depends only on the ellipticity constant of $a$, on $\|a\|_{C^{0,\alpha}(B_2)}$, on $d$ and on
$s,q^*$. In particular it does not depend on $u$ nor on $x_0$. Applying the De Giorgi-Nash estimate, we have
\begin{displaymath}
  \|u\|_{L^\infty(B_2(x_0))} \leq C \left(\int_{B_2(x_0)} u^2\right)^{1/2},
\end{displaymath}
where $C$ does not depend on $x_0$, for the same reasons as above. Hence, applying \eqref{eq:estimation_brezis} for
$s=2$, $u$ is bounded, which, by Liouville theorem (see for instance
\cite{Moser}), implies that $u$ is constant. 
\hfill $\diamondsuit$

\medskip

\begin{remark} In the proof of Proposition~\ref{prop:divergence-form}, we have used some of our results established in the case $q=2$
  in~\cite{milan}. We actually only made use of the uniqueness result (in order to prove that $\nabla u\equiv 0$
  in~\eqref{eq:zero-1} above), while we did establish~\cite{milan} existence and uniqueness of the solution
  (although not stated as such, the proofs of \cite{milan} and \cite{cpde-defauts} imply
  continuous dependency on the datum). If we allow ourselves to also use the existence result, then the above proof
  may be slightly simplified. One may then prove by continuation "only" that
  estimate~\eqref{eq:estimee-informelle-div} holds true (in particular, it in turn implies uniqueness), while the
  existence part is a consequence of a density argument: we approximate ~$f\in \left(L^q(\RR^d)\right)^{d}$ by a
  sequence~$f_n\in \left(L^2\cap L^q(\RR^d)\right)^{d}$; for each~$f_n$, we have a solution~$u_n$ the gradient of
  which is in~$\left(L^2(\RR^d)\right)^{d}$; using the estimate, the sequence $\nabla u_n$ is a Cauchy sequence
  in~$\left(L^q(\RR^d)\right)^{d}$; we may finally pass to the limit and obtain a solution for $f$ only
  in~$\left(L^q(\RR^d)\right)^{d}$. We chose to present the proof of Proposition~\ref{prop:divergence-form}
  because its pattern is more general and applies (see Section~\ref{sec:non-divergence}) to operators that are not in divergence form (to which our arguments of~\cite{milan} do not carry over). 
\end{remark}

\medskip

\begin{remark}\label{rk:L1}
It is well known that, even when $a=a^{per}\equiv 1$, estimate~\eqref{eq:estimee-informelle-div} is wrong for $q=1$
(in dimensions~$d>1$). For instance, for $d\geq 3$, let us define
\begin{displaymath}
  u(x) = \frac 1 {d|B_1|}\int_{B_1} \frac {(x-y)\cdot e} {|x-y|^{d}} dy,
\end{displaymath}
where $B_1$ is the unit ball of $\RR^d$, $|B_1|$ its volume, and $e\neq 0$ is a fixed vector in $\RR^d$. It is clear that
$-\Delta u = \div\left(\mathbf{1}_{B_1}e\right)$ in the sense of distributions. However, a simple computation shows
that $\nabla u(x) \approx  \frac{1}{d|x|^d}\left(e - d \frac{(x\cdot e)x}{|x|^2} \right),$ as $|x|\to +\infty$. Hence $\nabla u \not\in L^1\left(\RR^d\right)$. 
\end{remark}

\medskip

\begin{remark} We shall see in Section~\ref{sec:homogenization} below that a consequence of
  Proposition~\ref{prop:divergence-form} is the existence of a corrector in the adequate functional space, and, in
  turn, a quantitative theory of homogenization where the rates of convergence may be made precise, both for the
  Green functions associated to the divergence operators and for the solutions to the homogenized
  problems. Actually, the existence of a corrector in the adequate functional space conversely implies
  Proposition~\ref{prop:divergence-form}. By the arguments introduced for the periodic case in~\cite{AL1987},
  further made precise in~\cite{KLS} for various boundary conditions, and adapted in~\cite{josien} to the case of
  divergence operators with perturbed periodic coefficients satisfying~\eqref{eq:hyp1}, it is indeed possible to
  establish, from the existence of a suitable corrector, the approximation properties for the Green function
  $G(x,y)$. Then, the arguments of~\cite{AL1991}, using the representation formula for the solution of \eqref{eq:equation-informelle-div}, can be replicated to obtain the results of Proposition~\ref{prop:divergence-form}.
\end{remark}

\medskip

\begin{remark}\label{rk:systeme_2}
The statement and proof of Proposition~\ref{prop:divergence-form} concern \emph{equations}, but an analogous result
holds for \emph{systems} in divergence form. Indeed, all our arguments carry over to systems, including the central
result by Avellaneda and Lin of~\cite{AL1991} (based on the results of~\cite{AL1987}) on the continuity of
operators with periodic coefficients, and the uniqueness result that is a consequence of our arguments
of~\cite{milan} when $\nabla u\in \left(L^2(\RR^d)\right)^{d}$. The
only point above at which we have used the fact that we deal with an equation (actually, applying the Harnack
inequality) is when we prove that $\nabla w_{p,per}\in L^\infty$. But this is also implied by the results of
Avellaneda and Lin \cite{AL1987}. However, as stated in Remark~\ref{rk:systeme_1}, for systems, we do not have a
direct proof of the case $q<2$, and only can prove it by a duality argument. 
\end{remark}

\begin{remark}
  \label{rk:gloria}
 All what we need in the above proof of Proposition~\ref{prop:divergence-form} is that (i) the gradient $\nabla
 w_{per}$ of the periodic corrector is $L^\infty$ (and the latter fact is, in particular, true when $a_{per}$ itsef
 is H\"older continuous), (ii) the coefficient $a$ is uniformly continuous (in order to be able to use the local
 elliptic regularity result~\eqref{eq:20}), (iii) the result of Avellaneda and Lin concerning the
 continuity~\eqref{eq:estimee-informelle-div} in the case of a periodic coefficient. The recent results of
 \cite{gloria} allow to provide a functional analysis setting for non H\"older coefficients.
\end{remark}
\section{Estimate for operators in non-divergence form}
\label{sec:non-divergence}

%
%

The purpose of this section is to prove the result analogous to that of Proposition~\ref{prop:divergence-form} in the case of the equation, not in divergence form,
\begin{equation}
\label{eq:equation-informelle-non-div}
-a_{ij}\partial_{ij}u=f\quad \text{in } \RR^d.
\end{equation}
Note that because of the specific form of~\eqref{eq:equation-informelle-non-div}, we may assume, without loss of generality, that the matrix-valued coefficient $a$ in~\eqref{eq:equation-informelle-non-div} is symmetric.  

\medskip

Our result is:

\begin{proposition}
\label{prop:non-divergence-form}
Assume \eqref{eq:aper+tildea}-\eqref{eq:hyp1}. Fix $1<q<+\infty$. Then, for all $f\in L^q(\RR^d)$, there exists~$u\in L^1_{loc}(\RR^d)$ such that $D^2\,u\in L^q(\RR^d)$, solution to equation~\eqref{eq:equation-non-div} namely
\begin{equation}
  \tag{\ref{eq:equation-informelle-non-div}}
-a_{ij}\partial_{ij}u=f\quad \text{in } \RR^d.
\end{equation}
Such a solution is  unique up to the addition of an affine function. In addition, 
there exists a constant $C_q$, independent on $f$ and $u$, and only depending  on $q$, $d$ and the coefficient $a$, such that   $u$ satisfies  \begin{equation}
  \label{eq:estimee-informelle-non-div}
\left\|D^2\,u\right\|_{\left(L^q(\RR^d)\right)^{d\times d}}\,\leq C_q\, \left\|f\right\|_{L^q(\RR^d)}. 
\end{equation}
\end{proposition}



\medskip

Proposition~\ref{prop:non-divergence-form} will be used in the next section (and also in~\cite{BLL-2017-2}) to
proceed with the homogenization of the equation
 \begin{equation}
  \label{eq:equation-non-div}
  -a_{ij}(x/\varepsilon)\,\partial_{ij} u^\varepsilon=f.
\end{equation}
The correctors associated to~\eqref{eq:equation-non-div} may be put to zero, since they are solutions to
$-a_{ij}\partial_{ij} (p.x+w_p(x))=0$. Because it is not in divergence form, the homogenization
of~\eqref{eq:equation-non-div} (or the precise understanding of the  behavior of its solution~$u^\varepsilon$ for $\varepsilon$ small) however requires to understand the adjoint problem defining the invariant measure associated to~\eqref{eq:equation-non-div}. The latter reads as
\begin{equation}
\label{eq:equation-adjoint}
-\partial_{ij} (a_{ij}\,m)=0, \quad \text{in } \RR^d,
\end{equation}
or equivalently, 
\begin{equation}
\label{eq:equation-adjoint-tilde}
-\partial_{ij} (a_{ij}\,\tilde m)=\partial_{ij} (\tilde a_{ij}\,m_{per}),\quad \text{in } \RR^d,
\end{equation}
 decomposing, in the same spirit as we decomposed the corrector earlier, the measure as~$m=m_{per}+\tilde m$ ,
 where $-\partial_{ij} (a^{per}_{ij}\,m_{per})=0$. The existence and uniqueness of $m_{per}$, under the constraints
 $m_{per}\geq 0$ and $\left\langle m_{per}\right\rangle = 1$, is proved in \cite{blp}. The existence of $\tilde m$ in the suitable functional space is readily related to Proposition~\ref{prop:non-divergence-form}. We will see the details in Section~\ref{sec:homogenization}.
 
\subsection{From the non-divergence form to the divergence form}
\label{ssec:non-div-to-div}

To start with and as a preparatory work (both for the proof of Proposition~\ref{prop:non-divergence-form} and the homogenization of equation~\eqref{eq:equation-non-div} in Section~\ref{sec:homogenization}), we recall here, for convenience of the reader, a classical algebraic manipulation (see e.g.~\cite{blp}) that transforms  an equation in non divergence form to an equation in divergence form provided an invariant measure (that is, a solution to the adjoint equation) exists and enjoys suitable properties. 
We perform the transformation here in full generality and abstractly, in the case of the general equation  
\begin{equation}
\label{eq:general-equation}
-a_{ij}\partial_{ij} u + b_i \partial_{i}u  = f,
\end{equation}
where  $a_{ij}$, $1\leq i,j\leq d$, $b_i$, $1\leq i\leq d$, are general coefficients. We will actually use
the transformation at several distinct stages of our work in the present article, and also in our forthcoming
article~\cite{BLL-2017-2}. The coefficients $a_{ij}$, $b_i$ will either be periodic, or include the local
perturbation. They will either be at scale one (meaning $a_{ij}(x)$), or be rescaled by $\varepsilon$ as in
$a_{ij}(x/\varepsilon)$, etc. The first-order coefficients $b_i$, $1\leq i\leq d$, will identically vanish (as in
the case here), or not (in~\cite{BLL-2017-2}). 

\medskip
%

Consider~\eqref{eq:general-equation}, posed on a (not necessarily) bounded domain $\Omega$ and supplied with some suitable boundary conditions (or conditions at infinity) we do not make precise in this formal generic argument. Assume that there exists a positive solution~$m$, actually bounded away from zero, $\inf m>0$, to the adjoint equation 
\begin{equation}
\label{eq:general-equation-adjoint}
-\partial_{i} (\partial_j a_{ij}\,m + b_i \,m)=0, 
\end{equation}
on the same domain, with boundary conditions that we do not make precise either, and suitably normalized. Multiplying~\eqref{eq:general-equation} by $m$, we obtain, without even using the specific properties of~$m$, that
\begin{equation}
\label{eq:general-equation-m}
  -\div\left(\overline a \nabla u\right) + \left(\overline b + \div\overline a\right)\, .\,\nabla u+
  = \overline f,
\end{equation}
where 
\begin{equation}
\label{eq:m-coeff}
\overline a = m\, a,\quad\overline b = m\, b,\quad\overline f = m\,f.
\end{equation}
Precisely because of~\eqref{eq:general-equation-adjoint},  $\overline b +
\div\overline a$ is divergence-free
\begin{equation}
\label{eq:div-free}
  \div( \overline b +
\div\overline a) = 0.
\end{equation}
Equation~\eqref{eq:div-free} (again formally) implies the existence of a skew-symmetric matrix $\mathcal B$ such
that 
\begin{equation}
\label{eq:B}
   \overline b +
\div\overline a = \div \mathcal B.
\end{equation}
In the particular case of dimension $d=3$, this is equivalent to the existence of a
vector field $B=(B_1,B_2,B_3)$ such that 
\begin{equation}
\label{eq:B-3d}
   \overline b +
\div\overline a = \curl B,
\end{equation}
where $\mathcal B$ and $B$ are related by
\begin{equation}
\label{eq:mathcal-B}
  \mathcal B =
  \begin{pmatrix}
    0 & -B_3 & B_2 \\ B_3 & 0 & -B_1 \\ -B_2 & B_1 & 0
  \end{pmatrix}.
\end{equation}
Using $\mathcal B$ and relation \eqref{eq:B}, it is then immediate to observe that 
\begin{equation}
   \div\left[u\,(\overline b +
\div\overline a )\right] =\div(u\,\div \mathcal B) = \div(\mathcal B\, \nabla u),
\end{equation}
and thus~\eqref{eq:general-equation-m} reads as the equation \emph{in divergence form}
\begin{equation}\label{eq:8}
  -\div\left(\mathcal A \nabla u\right) = \overline f.
\end{equation}
with 
\begin{equation}
\label{eq:mathcal-A}
\mathcal A= \overline a- \mathcal B
\end{equation}
Since $\inf m>0$,  $\overline a$ is elliptic. Moreover, $\mathcal B$ is skew-symmetric, hence the matrix $\mathcal A$ is elliptic.

\medskip

As mentioned earlier, we will make the above transformation explicit, and justify it, in each specific instance we need. The first of  these instances, and actually a very classical and well known one, is a simple periodic setting. Consider~\eqref{eq:general-equation} posed on the entire space $\RR^d$  for periodic second-order coefficients $a_{ij}=a^{per}_{ij}$, $1\leq i,j\leq d$, with period the unit cell of the periodic lattice $\ZZ^d$, $b_i\equiv 0$, $1\leq i\leq d$, and $c\equiv 0$. The adjoint equation~\eqref{eq:general-equation-adjoint} to be considered is posed also on the entire space, for \emph{periodic} solutions,   and reads as~
\begin{equation}
\label{eq:inv-adj-per}
-\partial_{ij}(a^{per}_{ij}\, m_{per})=0.
\end{equation}
It is established, e.g. in~\cite{blp}, that there exists a unique nonnegative periodic solution $m_{per}$ that is
normalized, regular and is indeed bounded away from zero. Performing the above manipulations, we note that
$\overline b\equiv 0$ and $\div a^{per} $ is of zero mean in~\eqref{eq:B}. Thus, the matrix~$\mathcal B=\mathcal
B^{per}$ may be assumed periodic, and we write the equation originally considered in the divergence form 
\begin{equation}
\label{eq:forme-div-periodique}
-\div\left(\mathcal A^{per} \nabla u\right)=m_{per} f.
\end{equation}
We notice that, because $b_i\equiv 0$, $1\leq i\leq d$,  here, $\div\mathcal B^{per}=\div(m_{per}\,a^{per})$, and therefore
\begin{equation}
\label{eq:div-A-nulle}
\div\mathcal A^{per}=\div(m_{per}\,a^{per})-\div\mathcal B^{per}=0,
\end{equation}
a property we shall use in the next section.

\subsection{Proof of Proposition~\ref{prop:non-divergence-form}}
\label{ssec:proof-non-divergence-form}

The proof of Proposition~\ref{prop:non-divergence-form} essentially follows the same pattern as that
of~Proposition~\ref{prop:divergence-form}. We again argue by continuation, (this time for all  $1<q<+\infty$ since,
in this case, the exponent $q=2$ does not play any specific role), and show that the interval defined
by~\eqref{eq:interval}  is again the entire interval~$[0,1]$. Of course, this time Property ${\mathcal P}$ is
based on the statements of Proposition~\ref{prop:non-divergence-form}  and not those of Proposition~\ref{prop:divergence-form}  any longer. 

\medskip

The fact that $0\in {\mathcal I}$ is a consequence of the results of~\cite[Theorem B]{AL1991}, precisely because of
the algebraic manipulations we recalled above, which allow to rewrite the equation under the conservative
form~\eqref{eq:forme-div-periodique}, with the specific property~\eqref{eq:div-A-nulle}. The local integrability
$u\in L^1_{loc}(\RR^d)$ is like in Section~\ref{sec:divergence} obtained by elliptic regularity.

\medskip

Next, we show that ${\mathcal I}$ is open (relatively to the interval~$[0,1]$). In order to do so, we proceed
exactly as in the proof of Proposition~\ref{prop:divergence-form}, writing the equation $-(a_t+\varepsilon\tilde
a)_{ij} \partial_{ij} u = f$ as
\begin{equation}\label{eq:2}
  D^2u = \phi_t\left(f +\varepsilon\tilde a_{ij} \partial_{ij} u\right),
\end{equation}
where $\phi_t$ is the application $f\mapsto D^2u,$ where $u$ is the solution to
\begin{displaymath}
  -\left(a_t\right)_{ij} \partial_{ij} u = f.
\end{displaymath}
Here again, the map appearing in \eqref{eq:2} is proved to be a contraction for $\varepsilon>0$ sufficiently small,
thereby showing existence and uniqueness of the solution, together with the continuity estimate. 


\medskip

Regarding the closeness of ${\mathcal I}$, the heart of the matter is,  similarly to the case of an operator in divergence form, to show that if we have ~$f^n\in \left(L^{q}(\RR^d)\right)^{d}$ and $u^n$ with~$D^2 u_n\in \left(L^{q}(\RR^d)\right)^{d\times d}$, such that 
\begin{equation}
\label{eq:contradic1-bis-nondiv}
-\,(a_t)_{ij}\partial_{ij} u^n= \,f^n\quad \text{in } \RR^d,\end{equation}
\begin{equation}
\label{eq:contradic2-nondiv}
\left\|f^n\right\|_{\left(L^q(\RR^d)\right)^{d}}\buildrel n\longrightarrow +\infty\over \longrightarrow 0,
\end{equation}
\begin{equation}
\label{eq:contradic3-nondiv}
\left\|D^2 u^{n}\right\|_{\left(L^q(\RR^d)\right)^{d\times d}}=1,\quad\hbox{\rm for all}\,n\in \NN.
\end{equation}
then we reach a contradiction. To this end, we first prove, using the same argument as in the proof of
Proposition~\ref{prop:divergence-form} and the result by Avellaneda and Lin \cite{AL1991} on the operator with periodic coefficient (this time in non divergence form), that
\begin{equation}
\label{eq:cc-1-nondiv}
\exists\, \eta>0, \quad \exists\, 0<R<+\infty,\quad \forall \,n\in\NN,\quad \left\|D^2 u^{n}\right\|_{\left(L^q(B_R)\right)^{d\times d}}\geq \eta>0.
\end{equation}
Because of the bound~\eqref{eq:contradic3-nondiv}, we may claim that, up to an extraction, $D^2 u^n$ weakly
converges in $\left(L^q(\RR^d)\right)^{d\times d}$, to some $D^2 u$. Passing to the limit in the sense of
distribution implies that $u$ is a solution to $-\,(a_t)_{ij}\partial_{ij} u= 0$. Hence,
\begin{displaymath}
  -\,(a_t)_{ij}\partial_{ij}\left(u^n - u\right) = f^n.
\end{displaymath}
The Poincar\'e-Wirtinger inequality and~\eqref{eq:contradic3-nondiv} imply that, up to the addition of an affine function, $u^n$ is bounded in
$W^{2,q}(B_R)$.  Applying Rellich Theorem, we know that, up to extracting a subsequence, $u^n$ converges strongly
in $L^q(B_R)$, for any $R>0$. Elliptic regularity results \cite[Theorem 9.11]{GT} then imply
\begin{equation}\label{eq:4}
  \left\|u^n-u\right\|_{W^{2,q}(B_R)} \leq C(R)\left(\left\|f^n\right\|_{L^q(B_{R+1})} +\left\|u^n-u\right\|_{L^{q}(B_{R+1)}}  \right).
\end{equation}
Here, the constant $C(R)$ depends on $R$, on the ellipticity constant of $a$ and of its $C^{0,\alpha}_{\rm unif}(\RR^d)$ norm, but
not on $f^n,u^n,u$.
The right-hand side of this inequality tends to $0$ as $n\to+\infty$, hence we have strong convergence of $u^n$ to
$u$ in $W^{2,q}(B_R)$. In particular, \eqref{eq:cc-1-nondiv} implies that $u\not\equiv 0$. Concluding the proof amounts to reaching a contradiction with $-\,(a_t)_{ij}\partial_{ij} u= 0$. This requires a significantly different proof
from the case of operators in divergence form, because here we cannot
bootstrap some $L^2$ integrability and use coerciveness to conclude. The proof, in the present case, relies on the
maximum principle. 

\medskip

We first give the end of the proof assuming $d\geq 3$. We will explain below how to adapt it to the case $d=2$.

First, we assume $q<d/2$. In such a case, we claim that
\begin{equation}
  \label{eq:13}
\forall\, n\in\NN \text{ such that } n < \frac d {2q}, \quad D^2 u \in \left(L^{s_n}(\RR^d)\right)^{d\times d}, \text{
  where } \frac 1 {s_n} = \frac 1 q - \frac {2n} d. 
\end{equation}
This is proved by induction on $n$. The case $n=0$ is true by assumption. If we assume that \eqref{eq:13} is true
for $n-1$, with $n<d/(2q)$, the
Gagliardo-Nirenberg-Sobolev inequality \cite[Section 5.6.1]{evans} and the fact that
$D^2u\in \left(L^{s_{n-1}}\left(\RR^d\right)\right)^{d\times d}$ imply that, up to the addition of an affine function, $\nabla u \in
\left(L^{s_{n-1}^*}(\RR^d)\right)^d$ and $u\in L^{s_{n-1}^{**}}(\RR^d)$, where
\begin{displaymath}
  \frac 1 {s_{n-1}^*} = \frac 1 {s_{n-1}} - \frac 1 d, \qquad \frac 1 {s_{n-1}^{**}} = \frac 1 {s_{n-1}^*} - \frac 1 d = \frac 1 {s_{n-1}} - \frac 2 d.
\end{displaymath}
In other words, $s_{n-1}^{**} = s_{n}$. We then apply \cite[Theorem 9.11]{GT} again (that is, inequality \eqref{eq:4} with $u^n=0$ and $f^n=0$), finding
\begin{displaymath}
  \int_{B_1(x_0)} |D^2u|^{s_n}\leq \|u\|_{W^{2,s_n}(B_1(x_0))}^{s_n} \leq C \|u\|_{L^{s_n}(B_2(x_0))}^{s_n} = C\int_{B_2(x_0)}|u|^{s_n},
\end{displaymath}
where $C$ does not depend on $u$ nor on the center $x_0$ of the balls. Summing up all these estimates for $x_0\in
\delta \ZZ^d$, with $\delta>0$ sufficiently small, we obtain \eqref{eq:13}. Next, we choose $n$ such that
\begin{displaymath}
  \frac d {2q} - 1 < n < \frac d {2q},
\end{displaymath}
which is always possible because $d/q>2$. Then, we have $s_n>d/2$. Hence, Morrey's Theorem \cite[Section
5.6.2]{evans} implies that $u\in C^{0,\alpha}_{\rm unif}(\RR^d)$. Since we also have $u\in L^{s_{n-1}^{**}}(\RR^d)$, we infer that
$u$ vanishes at infinity: for any $\delta>0$, we have,
for $R$ sufficiently large, $|u(x)|<\delta$ if $|x|>R$. Applying the maximum principle, we infer that $-\delta \leq u
\leq \delta$ in $\RR^d$. This being valid for any $\delta>0$, we find $u\equiv 0$,
reaching a final contradiction.

Second, we assume that $q\geq d/2$. We claim that
\begin{equation}
  \label{eq:14}
\forall\, n\in\NN \text{ such that } \frac 1 q + \frac n r  < 1, \quad D^2 u \in \left(L^{\sigma_n}(\RR^d)\right)^{d\times d}, \text{
  where } \frac 1 {\sigma_n} = \frac 1 q + \frac n r. 
\end{equation}
Here again, we prove this by induction: for $n=0$, we have $\sigma_0 = q$ and
assumption~\eqref{eq:contradic3-nondiv} implies $D^2u \in L^q(\RR^d)^{d\times d}$. Assuming
that \eqref{eq:14} holds for $n-1$, with $n< r - r/q$, we write the equation satisfied by $u$ as
\begin{displaymath}
  a_{ij}^{per}\partial_{ij} u = t\tilde a_{ij} \partial_{ij}u \in L^{s}(\RR^d), \quad \frac 1 {s} = \frac 1
  {\sigma_{n-1}} + \frac 1 r = \frac 1 {\sigma_n},
\end{displaymath}
since $\tilde a \in L^r(\RR^d)$ and $D^2 u \in L^{\sigma_{n-1}}(\RR^d)$. Applying the results of \cite{AL1991}, we
thus have $D^2u \in \left(L^{\sigma_n}(\RR^d)\right)^d$. This concludes the proof of \eqref{eq:14}. 

Next, we choose $n$ such that
\begin{displaymath}
   r\left(\frac 2 d - \frac 1 q\right) < n < r\left(1 - \frac 1 q\right).
\end{displaymath}
This is possible if $r\left(1 - \frac 1 q\right) - r\left(\frac 2 d - \frac 1 q\right)>1$, that is, $r> \frac
d{d-2}$. Since $\tilde a \in L^r\cap L^\infty$, we may in fact increase $r$ so that this condition is
fulfilled. For such a value of $n$, we have $\sigma_n < d/2,$ and we may therefore apply our argument of the case $q<d/2$.
Here again, we reach a
contradiction. 

\medskip

Let us now assume that $d=2$. We cannot, as we did above, assume that $q<d/2$. However, the proof of \eqref{eq:14}
is still valid. We apply this inequality for the largest possible value of $n$, that is, $n = \left\lfloor r -
  \frac r q \right\rfloor$, where $\lfloor \cdot \rfloor$ denotes the integer part. We thus have
\begin{displaymath}
  D^2u \in \left(L^\sigma(\RR^d)\right)^{d\times d}, \quad \sigma = \sigma(r,q) = \frac{q}{1 + \frac q r \left\lfloor r -
  \frac r q \right\rfloor}.
\end{displaymath}
As we already pointed out above, since our assumption is that $\tilde a\in L^r\cap L^\infty(\RR^d)$, we may
increase $r$ if we wish. Since $\displaystyle \frac q r \left\lfloor r -
  \frac r q \right\rfloor \to q-1$ as $r\to+\infty$, we infer that
\begin{equation}\label{eq:35}
  D^2u \in \left(L^\sigma(\RR^d)\right)^{d\times d}, \quad\forall \sigma \in ]1,q].
\end{equation}

Next, we note that, since the ambient dimension is $d=2$, the fact that the matrix $a$ is elliptic and symmetric
implies that there exists two positive constants $C_0$ and $C_1$ such
that, for any symmetric matrix $e=e_{ij}$,
\begin{displaymath}
  C_0 \left(a_{ij} e_{ij}\right)^2 \geq  e_{ij}^2 + C_1 \det(e),
\end{displaymath}
with summation over repeated indices. 
This inequality is easily proved by elementary considerations, and was used for instance in \cite{cordes}, and
stated in \cite[Equation (4)]{talenti}. We apply
it to $e = D^2u$, multiply by $\chi_R^2$, where $\chi_R$ is a smooth cut-off function such that $\chi_R = 1$ in $B_R$, $\chi_R = 0$ in
$B_{R+1}^C$, and $|\nabla \chi_R|\leq 2$. We integrate over $\RR^2$ and find
\begin{displaymath}
  0 = C_0 \int \left(a_{ij}\partial_{ij}u\right)^2\chi_R^2 \geq \int \left(\partial_{ij}u\right)^2 \chi_R^2 +
  C_1\int\left(\partial_{11}u\partial_{22}u - \left(\partial_{12}u\right)^2\right)\chi_R^2.
\end{displaymath}
We note that the integrand in the last term is equal to $\partial_1\left(\partial_1 u\partial_{22}u \right)\chi_R^2
-\partial_2\left(\partial_1 u \partial_{12}u\right)\chi_R^2$. Integrating by parts, we thus have
\begin{displaymath}
  0 \geq \int \left(\partial_{ij}u\right)^2 \chi_R^2 - 2C_1 \int \chi_R\, \partial_1 u\, \partial_{22} u\, \partial_1
  \chi_R + 2C_1 \int \chi_R\, \partial_1u\, \partial_{12} u\, \partial_2 \chi_R.
\end{displaymath}
Hence,
\begin{equation}\label{eq:34}
  \int \left(\partial_{ij}u\right)^2 \chi_R^2 \leq 2C_1 \int |\nabla \chi_R||\nabla u| \chi_R |D^2 u|.
\end{equation}
If $q\geq 4/3$, \eqref{eq:35} implies $D^2 u\in L^{4/3}$, and, by the Gagliardo-Nirenberg-Sobolev inequality, $\nabla u \in
L^{4}$. As a consequence, $|\nabla u| |D^2 u|\in L^1$, and, letting $R\to +\infty$ in \eqref{eq:34}, we infer that
$D^2 u = 0$. If $q<4/3$, then, by the Gagliardo-Nirenberg-Sobolev inequality, $\nabla u \in L^{q^*}$, where $\frac 1
{q^*} = \frac 1 q - \frac 12$. In particular, the conjugate exponent of $q^*$, denoted by $(q^*)'$, satisfies $4/3
< (q^*)' < 2$. Thus, $q < (q^*)'< 2$. Hence, successively applying H\"older inequality and the interpolation inequality to \eqref{eq:34}, we infer
\begin{multline}
  \int |D^2 u|^2 \chi_R^2 \leq 2C_1 \bigl\| |\nabla \chi_R|\, |\nabla u|\bigr\|_{L^{q^*}(\RR^d)} \left\|\chi_R |D^2 u|
  \right\|_{L^{(q^*)'}(\RR^d)}, \\
\leq 2C_1 \bigl\||\nabla \chi_R|\, |\nabla u|\bigr\|_{L^{q^*}(\RR^d)} \left\|\chi_R |D^2 u|
  \right\|_{L^2(\RR^d)}^{\beta} \left\|\chi_R |D^2 u|
  \right\|_{L^q(\RR^d)}^{1-\beta},  
\end{multline}
for $\displaystyle \beta = \frac{4-3q}{2-q}.$
Thus,
\begin{displaymath}
  \left\|\chi_R |D^2 u|
  \right\|_{L^2(\RR^d)}^{2-\beta} \leq 2C_1 \bigl\||\nabla \chi_R|\, |\nabla u|\bigr\|_{L^{q^*}(\RR^d)}  \left\|\chi_R |D^2 u|
  \right\|_{L^q(\RR^d)}^{1-\beta} \mathop{\longrightarrow}_{R\to+\infty} 0,
\end{displaymath}
since $\nabla u \in L^{q^*}$ and $D^2u\in L^q$. Recalling that $\beta\in]0,1[$, thus $2-\beta>0$, we obtain $D^2u
= 0$.

Once again, we have reached a contradiction. This shows that  ${\mathcal I}$ is closed. As it is also open and non empty, it is
equal to~$[0,1]$ and this concludes the proof of Proposition~\ref{prop:non-divergence-form}. \hfill $\diamondsuit$





 
\medskip

\begin{remark}\label{rk:systeme_3}
In sharp contrast to the case of operators in divergence form, Proposition~\ref{prop:non-divergence-form} and its proof as presented above cannot be  extended to the case of \emph{systems}. In particular, we have made use of the result by Avellaneda and Lin for non-divergence form operators, which is, to the best of our knowledge, specific to equations. 
In addition, even though \emph{some} systems satisfy the maximum principle, we do not see how to adapt our proof of
uniqueness to the generic case of systems. Note also that, besides the usefulness of
Proposition~\ref{prop:non-divergence-form} on its own, the specific use we will make of that proposition in
homogenization theory is exposed in  Section~\ref{sec:homogenization}. The treatment of non-divergence form operators there will require the use of the invariant measure associated to their adjoint, a concept we do not even know how to define for systems.
\end{remark}

\section{Application to homogenization}
\label{sec:homogenization}

\subsection{Divergence form}
\label{ssec:divergence}

We return to the corrector equation ~\eqref{eq:correcteur}, namely 
$$
  -\,\hbox{\rm div}\left(a\,(p+\nabla w_{p})\right)=0 \quad \text{in } \RR^d
$$
which we write under the form~\eqref{eq:correcteur-tilde}:
$$
  -\,\hbox{\rm div}\left(a\,\nabla \tilde w_{p}\right)= \,\hbox{\rm div}\left(\tilde a\,(p+\nabla w_{p,per})\right)\quad \text{in } \RR^d.
$$
Since $w_{p,per}$ is the periodic corrector, that is the solution to~\eqref{eq:correcteur-per}
$$  -\,\hbox{\rm div}\left(a^{per}(x)\,(p+\nabla w_{p,per}(x))\right)=0,$$
with the coefficient $a^{per}$ satisfying assumptions~\eqref{eq:hyp1}, we have, as pointed out at the beginning of
Section~\ref{sec:divergence}, $\nabla w_{p,per}\in \left(L^\infty(\RR^d)\right)^d$. 
We insert this information in the right-hand side of~\eqref{eq:correcteur-tilde}, and may therefore conclude, using Proposition~\ref{prop:divergence-form}
 for the specific exponent $q=r$ that there exists a function $\tilde w_{p}$, uniquely defined up to the addition of a constant, that solves~\eqref{eq:correcteur-tilde}, with  $\nabla \tilde w_{p}\in \left(L^r(\RR^d)\right)^d$ and with, considering~\eqref{eq:estimee-informelle-div},
 $$ \left\|\nabla \tilde w_p\right\|_{\left(L^r(\RR^d)\right)^{d}}\,\leq C_r\, \left(|p|+\left\|\nabla
     w_{p,per}\right\|_{\left(L^\infty(\RR^d)\right)^{d}}\right)\, \left\|\tilde
   a\right\|_{\left(L^r(\RR^d)\right)^{d}}.$$
Setting $w_p=w_{p,per}+\tilde w_p$, returning to equation~\eqref{eq:correcteur} and using the regularity~\eqref{eq:hyp1}, we also have that $\nabla \tilde w_p\in \left(L^\infty(\RR^d)\right)^d$.
We have therefore provided an alternative proof of our main results in Theorem 4.1 of~\cite{cpde-defauts}. The arguments of~\cite{josien} then allow to prove quantitative
 homogenization results. Recall however Remark~\ref{rk:notL1}: the above argument does not cover the case $r=1$, since in
Proposition~\ref{prop:divergence-form}, $q=1$ is excluded. 

\medskip

To end this section, let us mention that, using the above computation, if $G$ is the Green function associated to
\eqref{eq:general-equation}, and if $\mathcal G$ is the Green function
associated to \eqref{eq:8}, we have
\begin{displaymath}
  \mathcal G (x,y) = m(y) G(x,y). 
\end{displaymath}
Moreover, the measure $m$ satisfies $\gamma \leq m(x)\leq \frac 1 \gamma$, for some $\gamma>0$, and $m\in
C^{0,\alpha}_{\rm unif}(\RR^d)$. Hence, all the estimates which are valid for the Green function $G$ give estimates on
$\mathcal G$. The same conclusion holds for $\nabla_x \mathcal G$. On the other hand, estimates on $\nabla_y
\mathcal G$ can only be proved if $m\in W^{1,\infty}$, which is in general not the case. 

\subsection{Non-divergence form}
\label{ssec:non-divrgence}

We discuss here homogenization for the equation in non-divergence form~\eqref{eq:equation-non-div}
 $$
   -a_{ij}(x/\varepsilon)\,\partial_{ij} u^\varepsilon=f.
$$


In order to deal with this problem, we may apply two different strategies, both relying on the central estimate of
Proposition~\ref{prop:non-divergence-form}: 
\begin{itemize}
\item The first one consists in using this estimate to derive a bound on the
distance between the solution $u^\varepsilon$ and the solution corresponding to the case $\tilde a = 0$. Then,
using the results of \cite{AL1987,AL1989}, valid only in the periodic case, we may prove equivalent convergence
estimates in the present case. 

\item The second one, which is the one we chose to apply below, consists in using
  Proposition~\ref{prop:non-divergence-form} to prove the existence of a stationary measure. Then, multiplying the
  equation by this measure $m$, the calculations performed in Subsection~\ref{ssec:non-div-to-div} allow to write the
  above equation in divergence form. We may therefore use the same method as in
  Subsection~\ref{ssec:divergence}. This strategy seems more intrinsic and more easily adaptable to different
  situations. We follow it here.
\end{itemize}

We claim there exists an invariant measure associated to this equation. The precise result is the following.
\begin{proposition}
\label{prop:invariant-measure}
Assume the coefficient $a$ in~\eqref{eq:equation-non-div} satisfies~\eqref{eq:hyp1} with $r>1$. There exists an invariant measure associated to this equation, that is, by definition, a unique function~$m$ solution to~\eqref{eq:equation-adjoint}:
\begin{equation}
  \label{eq:mesure-invariante}
-\partial_{ij} (a_{ij}\,m)=0, \quad \text{in } \RR^d,  
\end{equation}
which writes $m=m_{per}+\tilde m$ where  $m_{per}$ is the periodic invariant measure $m_{per}$  (defined in
Section~\ref{ssec:non-div-to-div} above), and~$\tilde m$ is defined as the unique solution in $L^r(\RR^d)$ to~\eqref{eq:equation-adjoint-tilde}
\begin{equation}
  \label{eq:mesure-invariante-tilde}
-\partial_{ij} (a_{ij}\,\tilde m)=\partial_{ij} (\tilde a_{ij}\,m_{per}),\quad \text{in } \RR^d.  
\end{equation}
  This function $m$ is H\"older continuous, positive, bounded away from zero. 
  \end{proposition}
  \begin{remark}
In the case $r=1$, the above result still holds, but we only have $\tilde m \in L^q(\RR^d)$, for any $q>1$.    
  \end{remark}

The existence and uniqueness of $m$ as stated in Proposition~\ref{prop:invariant-measure} is an immediate
consequence of the following corollary of Proposition~\ref{prop:non-divergence-form}, the proof of which is
postponed until the end of this Section. 

\begin{corollary}[of Proposition~\ref{prop:non-divergence-form}]
\label{cor:adjoint-non-div}
Assume \eqref{eq:aper+tildea}-\eqref{eq:hyp1}. For all~$1<q<+\infty$ and $f\in\left(L^q(\RR^d)\right)^{d\times d}$, there exists a unique $u\in L^q(\RR^d)$ solution, at least in the sense of distributions,  to 
\begin{equation}
\label{eq:cor1}
-\,\partial_{ij} \left ( \,a_{ij}\,u\right)=\partial_{ij}\,f_{ij} \quad \text{in } \RR^d.
\end{equation}
This function $u$ satisfies 
\begin{equation}
\left\|u\right\|_{L^q(\RR^d)}\leq {\mathcal C}_{q}\,\left\| f\right\|_{\left(L^q(\RR^d)\right)^{d\times d}},
\end{equation}
for a constant ${\mathcal C}_q$ independent of $f$. 
\end{corollary}

\noindent{\bf Proof of Proposition~\ref{prop:invariant-measure}.} Applying Corollary~\ref{cor:adjoint-non-div}, we
know that there exists a solution $\tilde m$ to \eqref{eq:mesure-invariante-tilde}, with $\tilde m\in
L^r(\RR^d)$. In the case $r=1$, we have $\tilde m \in L^q(\RR^d)$, for any $q>1$. The measure $m_{per}$ solution to
$\partial_{ij}\left(a_{ij}^{per} m\right) = 0$ is already known to exist (see \cite{engquist-souga}), to be H\"older
continuous thanks to standard elliptic regularity results, and to be bounded away from $0$. The measure $m = m_{per} +
\tilde m$ is a solution to \eqref{eq:mesure-invariante}. 
We prove now that $m$ is positive. For this purpose, we first point out that $\tilde m$ is uniformly H\"older continuous,
thanks to the results of \cite{bogachev,bogachev-livre}. Since it is in $L^r(\RR^d)$, we know that
\begin{displaymath}
  \left\| \tilde m\right\|_{L^\infty(B_R^c)} \mathop{\longrightarrow}^{R\to +\infty} 0.
\end{displaymath}
Hence, for $R$ sufficiently large, we have 
\begin{equation}\label{eq:3}
  \forall\, x\in B_R^c, \quad m(x) \geq \frac12 \inf m_{per} >0.
\end{equation}
Applying the maximum principle on $B_R$, we infer that $m\geq 0$ in the whole space $\RR^d$. Next, we
apply the Harnack inequality \cite{bogachev}, which implies that $m$ is bounded away from $0$. This concludes the
proof of Proposition~\ref{prop:invariant-measure}.\hfill $\diamondsuit$

\medskip

Next, we rescale $m$, considering $m_\varepsilon(x)=m(x/\varepsilon)$ and
multiply~\eqref{eq:equation-non-div} by $m_\varepsilon$. The standard manipulations
(\eqref{eq:general-equation-m} through~\eqref{eq:mathcal-A}) recalled in Section~\ref{ssec:non-div-to-div} yield 
\begin{equation}
\label{eq:homog-div-1}
  -\div\left(\mathcal A_\varepsilon \nabla u^\varepsilon \right) = m_\varepsilon  f,
\end{equation}
with the elliptic matrix valued coefficient $\mathcal A_\varepsilon(x)=\mathcal A(x/\varepsilon)$, 
\begin{equation}
\label{eq:homog-mathcal-1}
\mathcal A = m\, a - \mathcal B
\end{equation}
and the skew-symmetric matrix-valued coefficient~$\mathcal B$ defined by~\eqref{eq:mathcal-B}.
In the specific case considered, where $m=m_{per}+\tilde m$, $\mathcal B$ is defined as the sum $\mathcal
B=\mathcal B^{per}+\tilde{\mathcal B}$, where the periodic part $\mathcal B^{per}$ is obtained solving the periodic
equation $\div \mathcal B^{per}=\div(m_{per}\,a^{per})$ (the right-hand side being divergence-free because of~\eqref{eq:inv-adj-per}, we recall) and where 
\begin{equation}
\label{eq:curlB}
\div \tilde{\mathcal B}=\div\left(\tilde m\,a^{per}+m_{per}\,\tilde a + \tilde m \, \tilde a\right).
\end{equation}
The latter equation (which also has a divergence-free right-hand side because of \eqref{eq:equation-adjoint-tilde},
that is, \eqref{eq:mesure-invariante-tilde}) admits a skew-symmetric solution $\tilde{\mathcal B}\in \left(L^r(\RR^d)\right)^{d\times d}$ which is unique up to the
addition of a constant. This is an application of the Calder\'on-Zygmund operator theory. Indeed, we introduce the solution to 
\begin{equation}
\label{eq:DeltaB}
-\Delta \tilde{\mathcal B}_{ij}=\partial_{jk}\Bigl(\tilde m a^{per}_{ik}+\left(m_{per}+\tilde m\right)\tilde
a_{ik}\Bigr) - \partial_{ik}\left(\tilde m a^{per}_{jk}+\left( m_{per}+\tilde m\right)\tilde a_{jk}\right),
\end{equation}
which is known to exist thanks to \cite{meyer,stein}, with the additional property that
the corresponding operator is continuous from $L^r$ to $L^r$:
\begin{multline}\label{eq:5}
  \left\|\tilde{\mathcal B}_{ij}\right\|_{L^r(\RR^d)} \leq C\sum_{k=1}^d \left(\left\|\tilde m
      a_{jk}^{per}+\left( m_{per}+\tilde m\right)\tilde a_{jk}\right\|_{L^r(\RR^d)}\right.\\ \left. + \left\|\tilde m a_{ik}^{per}+\left( m_{per}+\tilde m\right)\tilde a_{ik}\right\|_{L^r(\RR^d)}\right),
\end{multline}
where $C$ is a universal constant. Now, using that $\div(\div(ma)) = 0,$ a simple computation gives
\begin{multline*}
  -\partial_i \Delta \tilde{\mathcal B}_{ij} = -\partial_{iik}\left(\tilde m a^{per}_{jk}+\left( m_{per}+\tilde
      m\right)\tilde a_{jk}\right) \\= -\partial_k \Delta \left(\tilde m a^{per}_{jk}+\left( m_{per}+\tilde m\right)\tilde a_{jk}\right),
\end{multline*}
hence the distribution $T = \div \tilde{\mathcal B}- \div\left(\tilde m\,a^{per}+m_{per}\,\tilde a + \tilde m
  \, \tilde a\right)$ is
harmonic. Since, according to \eqref{eq:5} and the fact that $\tilde a,\tilde m\in L^r(\RR^d)$, $T\in
W^{-1,r'}(\RR^d)$, we necessarily have $T=0$, hence $\tilde{\mathcal B}$ satisfies \eqref{eq:curlB}. Finally, we
point out that the regularity assumed on $a_{per}$ and $\tilde a$ implies that $m_{per}$ and $\tilde m$ are both H\"older
continuous, and consequently that $\mathcal A = ma -\mathcal B$ satisfy the assumptions~\eqref{eq:hyp1}. 
On the other hand, the right-hand side $m_\varepsilon\,f$ of~\eqref{eq:homog-div-1} strongly converges (to~$f$) in $H^{-1}(\RR^d)$ as $\varepsilon$ vanishes. We may therefore apply the results of~\cite{cpde-defauts,josien} to~\eqref{eq:homog-div-1} and obtain the homogenized limit, with actual rates of convergence. 

\medskip

We conclude this section with the proof of Corollary~\ref{cor:adjoint-non-div}.

\medskip

\noindent {\bf Proof of Corollary~\ref{cor:adjoint-non-div}}
We fix $f\in\left(L^q(\RR^d)\right)^{d\times d}$, for some~$1<q<+\infty$, and denote by $q'$ the conjugate exponent of $q$, that is, $\displaystyle {1\over q}+{1\over {q'}}=1$. To any arbitrary function~$g\in L^{q'}(\RR^d)$, we may associate the unique (up to the addition of an affine function) solution~$v$, such that $D^2 v\in\left(L^{q'}(\RR^d)\right)^{d\times d}$, 
to $-\,a_{ij}\,\partial_{ij} v= \,g$. The map
\begin{displaymath}
  \begin{array}{rcl}
    L^{q'}(\RR^d) & \longrightarrow & \left(L^{q'}(\RR^d)\right)^{d\times d} \\
    g &\longmapsto & \partial_{ij}v
  \end{array}
\end{displaymath}
is linear continuous. We may therefore define by $g\longmapsto L_f(g):=\int_{\RR^d} f_{ij}\,.\,\partial_{ij} v$ a linear form
on $L^{q'}(\RR^d)$. Since, 
using the result of Proposition~\ref{prop:non-divergence-form} for $q'$,
\begin{eqnarray}
\label{eq:linear-form-1-nondiv}
\left | L_f(g)=\int_{\RR^d}  f_{ij}\,.\,\partial_{ij} v\right |& \leq &\,\left\| f\right\|_{\left(L^q(\RR^d)\right)^{d\times d}}\,
\left\| D^2 v\right\|_{\left(L^{q'}(\RR^d)\right)^d}\nonumber\\
&\leq& C_{q'}\,\left\| f\right\|_{\left(L^q(\RR^d)\right)^{d\times d}}\,
\left\| g\right\|_{L^{q'}(\RR^d)},
\end{eqnarray}
 this linear form is therefore a continuous map on $L^{q'}(\RR^d)$. Hence there exists some $u\in L^q(\RR^d)$ such that 
 $$L_f(g)=\int_{\RR^d}f_{ij}\,.\,\partial_{ij} v=\int_{\RR^d} g\,u,$$
and we read on estimate~\eqref{eq:linear-form-1-nondiv} that 
$$\left\|u\right\|_{L^q(\RR^d)}\leq C_{q'}\,\left\| f\right\|_{\left(L^q(\RR^d)\right)^{d\times d}}.
$$
There remains to show that $u$ satisfies $-\,a_{ij}\,\partial_{ij} u=\,\partial_{ij} f_{ij}$. To this end, we
consider the specific case where $v\in \mathcal{D}(\RR^d)$ (that is, $v$ is smooth and has compact support) and set $g=-a_{ij}\,\partial_{ij} v$.
Applying the above, we have  $\displaystyle \int_{\RR^d} f_{ij}\,.\,\partial_{ij} v=\int_{\RR^d} g\,u$. The left-hand side is the duality  product $\langle\,\partial_{ij}\,f_{ij},v\rangle_{{\mathcal D'}(\RR^d),\mathcal{D}(\RR^d)}$, while, by definition of~$g$, the right-hand side reads as 
$$ \int_{\RR^d} -a_{ij}\,\partial_{ij} v\, u=\langle\,-\,a_{ij}\,.\,\partial_{ij} u,v\rangle_{{\mathcal D'}(\RR^d),\mathcal{D}(\RR^d)}.$$
 Since this holds true for all~$v\in \mathcal{D}(\RR^d)$, this shows  $-\partial_{ij} \left ( \,a_{ij}\,u\right)=\partial_{ij}\,f_{ij}$ and concludes our proof.\hfill$\diamondsuit$



%

\section*{Acknowledgement} The work of the second author is partly
  supported by  ONR under Grant N00014-15-1-2777  and by EOARD, under Grant FA9550-17-1-0294. The authors wish to
  thank Marc Josien for pointing out the mistake of \cite{cpde-defauts} that is corrected in the Appendix below.
  
\appendix

\section{Appendix: Corrigendum to \cite{cpde-defauts}}
\label{sec:Appendix}


The present section aims at correcting some mistakes in \cite{cpde-defauts}. 

First, we recall Remark~\ref{rk:notL1} above that was pointing out that the case $r=1$ should not have been included in
Theorem~4.1 of \cite{cpde-defauts}.


Second, the purpose of this appendix is to correct \cite[Lemma 4.2]{cpde-defauts}. We have claimed there that the elementary  pointwise estimations known on the first and second gradient of the Green function of the Laplace operator can be generalized in similar estimates integrated locally for the operator~$-\hbox{\rm div}\, \left(a\,\nabla .\right)$ . 
Our precise statement, and some parts of the proofs are erroneous (however, this does not affect the
other results in \cite{cras-defauts,cpde-defauts}). In this appendix, we are going to use the notation of
\cite{cras-defauts,cpde-defauts}: $a_{per}$ is replaced by $a_0$, and $\tilde a$ by $b$. 

The correct statement of our results is as follows.

\begin{lemma}[Corrected version of {\cite[Lemma 4.2]{cpde-defauts}}]
\label{lem:1}
Assume that the coefficient~$a$ satisfies $a=a_{0}+b$, where $a_0$ denotes the (unperturbed) background, and $b$ the perturbation. Assume that $0<\underline{\mu}\leq a_{0}(x)+b(x)$, $0<\underline\mu \leq a_0(x)$, a.e.,  for some fixed constant $\underline{\mu}$, $a_{0}\in L^\infty(\RR^d)$, $b\in L^\infty(\RR^d)$. Consider the  Green function~$G$, solution to
\begin{equation}
-\hbox{\rm div}_x\, \left(a(x)\,\nabla_x G(x,y)\right)=\delta(x-y)
\end{equation}

\noindent (i) Then, for all $1\leq q\leq 2$, there exists a constant $C$ such that, for all~$R>0$ and all~$x\in \RR^d$, $G$ satisfies
\begin{equation}
\label{eq:grad-G}
\int_{B_{2R}(x)\backslash B_R(x)} \left|\nabla_yG(x,y)\right|^q\,dy\leq {\frac
  C {R^{d(q-1)-q}}},
\end{equation}
where $B_{2R}(x)\backslash B_R(x)=\left\{y, \ R\leq |x-y|\leq 2R\right\}$ denotes the annular region enclosed between the balls of radius~$R$ and $2R$.

\noindent (ii) Assume in addition that $a_0=a_{per}$ is periodic and H\"older continuous, and that $ b\in
L^r(\RR^d)\cap C^{0,\alpha}_{\rm unif}(\RR^d)$, for some $1\leq r<+\infty$, then $G$ satisfies 
\begin{equation}
\label{eq:grad-grad-G-new}
\forall\, q\in]1,+\infty[, \quad\exists\,  C>0, \quad \forall\, y\in\RR^d,\quad\int_{\{|x-y|>1\}}|\nabla_x\,\nabla_yG(x,y)|^qdx \leq C.
\end{equation}
\end{lemma}
Three comments are in order:
\begin{description}
\item[{[a]}] the estimate~\eqref{eq:grad-G} is correctly stated in~\cite{cpde-defauts} (as estimate (26) therein),
 but the proof there has a flaw. For clarity, we provide the \emph{entire}, corrected proof here. In the course of
 the proof of~\eqref{eq:grad-G} in~\cite{cpde-defauts}, it is indeed claimed  that the estimate 
 \begin{equation}
   \label{eq:estimee_grad_fausse}
\int_{B_{2R}\backslash B_R} \left|\nabla_xG(x,y)\right|^q\,dy\leq {\frac
  C{R^{d(q-1)-q}}}   
 \end{equation}
 holds true. 
  Note the gradient in $x$ and not in $y$ in the integrand. It is actually unclear that the latter estimate is correct, and we suspect it is not.
\item [{[b]}] it is claimed in~\cite{cpde-defauts}
that 
the estimate
\begin{equation}\label{eq:estimee_grad_grad_fausse}
  \int_{B_{2R}(x)\backslash B_R(x)}\left|\nabla_x\,\nabla_yG(x,y)\right|^q\,dy\leq {\frac C{R^{d(q-1)}}}
\end{equation}
holds; we are only able to establish this estimation as a consequence of a more precise, namely pointwise, estimation of $\nabla_x\,\nabla_yG(x,y)$ which is a consequence of arguments in both~\cite{cpde-defauts} and~\cite{josien} and provided the additional assumption $ b\in L^r(\RR^d)$, for some $1\leq r<+\infty$ holds.
\item [{[c]}] we therefore replace this estimate on annular regions by the estimate~\eqref{eq:grad-grad-G-new}, and
  show it is sufficient to conclude the proof of Lemma~4.2 in~\cite{cpde-defauts}, thereby checking there is no circular argument.
\end{description}

\noindent {\bf Proof:} 
(i) The proof of \eqref{eq:grad-G} is exactly that of \cite{cpde-defauts}: we first prove that, for $1\leq q \leq 2$, 
\begin{equation}
  \label{eq:echange-dx-dy}
  \int_{B_{2R}(y)\backslash B_R(y)} \left|\nabla_xG(x,y)\right|^q\,dx\leq {\frac
  C{R^{d(q-1)-q}}}\ .
  \end{equation}
Note that in \eqref{eq:grad-G} and \eqref{eq:echange-dx-dy}, the role played by $x$ and $y$ are reversed. We will
see below that \eqref{eq:echange-dx-dy} indeed implies \eqref{eq:grad-G}. 

\medskip

  We first proceed for dimensions~$d\geq 3.$  In the case $q=2$, we use the Caccioppoli inequality (see e.g.~\cite[Lemma~4.3]{cpde-defauts}), which implies
that
$$\int_{B_{R/2}(x_0)} |\nabla_x G(x,y)|^2 dx \leq \frac C {R^2}
\int_{B_{R}(x_0)} |G(x,y)|^2 dx,$$
for any $x_0\in \RR^d$ such that $y\notin B_R(x_0)$.
Next, we fix~$y$ and we cover $B_{2R}\setminus B_R = \left\{x, \ R<|x-y|<2R \right\}$ by balls $B_{R/2}(x_i)$, for  some 
points $x_i$ such that $5R/4<|x_i|<7R/4$, in such a
way that (i) a finite number of such $x_i$ is sufficient to cover the ring $B_{2R}\setminus B_R$ and that (ii) any point in $B_{2R}\setminus B_R$ belongs to at most  $K$ balls $B_{R/2}(x_i)$, for some  $K$ that is independent of the radius~$R$. This is easily seen to be possible.

The above estimate holds
for any couple of balls $(B_{R/2}(x_i),B_{R}(x_i))$. We sum all such estimates over the finite number of indices~$i$ and obtain
\begin{equation}
  \label{eq:caccioppoli_appliquee}
  \int_{B_{2R}\setminus B_R} |\nabla_x G(x,y)|^2 dx \leq
\frac{C\,K}{R^2} \int_{B_{11R/4}\setminus B_{R/4}} |G(x,y)|^2 dx.
\end{equation}
Since $d\geq 3$, using the classical pointwise estimate
\begin{equation}
  \label{eq:estimation_G}
\forall\, x,y\in \RR^d, \quad 0\leq G(x,y)\leq \frac C {|x-y|^{d-2}}.  
\end{equation}
(established in~\cite{GW,KN} and recalled in~\cite[estimate (28)]{cpde-defauts}), we get, \begin{equation}
\label{eq:anneau}
\int_{B_{2R}\setminus B_R} |\nabla_x G(x,y)|^2 dx \leq \frac {C\,K} {R^2}\int_{R/4}^{11R/4}
\frac{r^{d-1}}{r^{2d-4}}dr \leq \frac{C}{R^{d-2}}.
\end{equation}
This proves the case $q=2$. For $q<2$, we simply apply the H\"older
inequality and use~(\ref{eq:anneau}):
\begin{eqnarray}
\label{eq:holder0}
\int_{B_{2R}\setminus B_R} |\nabla_x G(x,y)|^q dx &\leq &
\left(\int_{B_{2R}\setminus B_R} |\nabla_x G(x,y)|^2 dx\right)^{q/2}
R^{d(1-q/2)} \nonumber \\ 
&\leq& C R^{-(d-2)q/2 +d - dq/2} = C R^{-dq+d+q} .
\end{eqnarray}
We thus have proved (\ref{eq:echange-dx-dy}) for $d\geq 3.$

\medskip

We next prove (\ref{eq:echange-dx-dy}) for $d=2$. For this purpose, we use the following
inequality, valid for any $\beta\in (0,2]$, and which expresses and quantifies the continuous embedding  of $L^{2,\infty}$ into $L^r$ for $r<2$ on bounded domains:
$$\forall\, f\in L^{2,\infty}(\Omega), \quad \int_\Omega |f|^{2-\beta}
\leq C_\beta\,|\Omega|^{\beta/2} \|f\|_{L^{2,\infty}(\Omega)}^{2-\beta},$$
where $C_\beta=4 \frac{1+2^{-\beta}}{(2^\beta-1)^{(2-\beta)/2}} $ is suitable. This estimate is proved for instance in the Appendix of
\cite{anantha}. We are going to apply it to $f=\nabla_x G$ and $\Omega =
B_{2R}\setminus B_R$. Since, in sharp contrast to the situation in dimensions~$d\geq 3$, $G(x,y)$ does not vanish when $|x-y|\longrightarrow +\infty$, we use the estimate
\begin{equation}
\label{eq:G_2d_2}
\left\|\nabla_x G(.,y)\right\|_{L^{2,\infty}} \,\leq\,C,
\end{equation}
to bound from above the right hand side.
We find:
$$\int_{B_{2R}\setminus B_R} |\nabla_x G|^{2-\beta} \leq C_\beta\,C R^{\beta}.$$
This implies (\ref{eq:echange-dx-dy}) for $q = 2-\beta \in [0,2).$ Finally,
in order to prove (\ref{eq:echange-dx-dy}) for $q=2$, we fix~$y$ and first point out
that, integrating the equation 
$-\div_x(a(x)\nabla_x G(x,y)) =
\delta_y(x)$
on the set $\{x, \ G(x,y)\geq s\} $ (which contains $y$) for some $s\in \RR$, that
\begin{equation}
\label{eq:ipp0}
1 = \int_{ G\geq s } -\div_x(a\nabla_x G)dx = -\int_{G=s} \left(a
  \nabla_x G\right)\cdot n_s,
\end{equation}
where $n_s$ denotes the outward normal to the set $\{x, \ G(x,y)\geq s\}.$ Note that, here, we have implicitly assumed that
the set $\{x, \ G(x,y)\geq s\} $ is Lipschitz-continuous, so that its outer normal is well defined and we can
integrate by parts. This may not be the case, given the regularity of $G$. However, using the co-aera formula (see \cite[Theorem 3.40]{afp}), 
it is simple to prove that, since $G\in
C^{0,\alpha}$ away from $x=y$, this set is Lipschitz-continuous for almost all $s\in \RR$. This is sufficient
for our purpose here.

Next, we
multiply the equation by $G$ and integrate on $\{m\leq G \leq M\}$ for
some $m\leq M$. This gives
\begin{eqnarray}
0 &=& \int_{M\geq G \geq m} -\div_x(a\nabla_xG)\,G \,dx\nonumber\\ 
&=& \int_{G=M} G \left(a\nabla_x G \right)\cdot n_M - \int_{G=m} G
\left(a\nabla_x G \right)\cdot n_m + \int_{M\geq G \geq m} \left(a\nabla_x
  G\right)\cdot \nabla_x G.\nonumber\\ 
\end{eqnarray}
Hence, using (\ref{eq:ipp0}), we have
$$\int_{M\geq G \geq m} \left(a\nabla_x
  G\right)\cdot \nabla_x G = M-m.$$
Next, we define, for $R>0$, $m_R = \inf \left\{G(x,y), \ x\in
  B_{2R}\setminus B_R\right\},$ and $M_R = \sup \left\{G(x,y), \ x\in
  B_{2R}\setminus B_R\right\}.$ We have 
$B_{2R}\setminus B_R \subset \left\{ m_R \leq G \leq M_R\right\}$.
Hence,
$$\int_{B_{2R}\setminus B_R} |\nabla_x G|^2 dx \leq C \int_{m_R \leq G
  \leq M_R}\left(a\nabla_x
  G\right)\cdot \nabla_x G = M_R - m_R.$$
We apply the estimate (30) of \cite{dolzmann}, namely here
$$\left\|G\right\|_{C^{0,\alpha}(B_r)}\leq C\,r^{-\alpha}\,\left\|\nabla G\right\|_{L^{2,\infty}(B_{2r})},$$
for all $r$ such that $B_{2r}\subset B_{2R}\setminus B_R$. 
In view of~(\ref{eq:G_2d_2}),  that estimate implies
that $M_R-m_R$ is bounded independently of $R$. This proves
(\ref{eq:echange-dx-dy}) in the case $q=2$.

\medskip

At this stage, we have proved (\ref{eq:echange-dx-dy}). We next point out
that, in all generality and for non necessarily symmetric matrix-valued coefficients~$a$, $H(x,y) = G(y,x)$ is the Green function of the operator $-\div(a^T \nabla
\cdot)$, where $a^T$ is the transpose matrix of $a$, which satisfies the
same assumptions as $a$. Hence, we may apply (\ref{eq:echange-dx-dy}) to
$H$, finding (\ref{eq:grad-G}).

\bigskip

\noindent(ii) We now turn to the proof of \eqref{eq:grad-grad-G-new}.

Here again, we first proceed with the case $d\geq 3$, and deal with $d=2$ separately. 

We note that $G$
satisfies $-\div_x\left(a\nabla_x G(\cdot,y)\right) = 0$ in the set $|x-y|>1/2.$ Hence, applying \cite[Theorem
8.32]{GT}, we have
\begin{multline}\label{eq:GT-theorem 8.32}
  \forall\, x_0\in \RR^d \text{ such that } |x_0-y|>1, \\ \|G(\cdot,y)\|_{C^{1,\alpha}(B_{1/4}(x_0))} \leq C \|G(\cdot,y)\|_{L^\infty(B_{1/2}(x_0))},
\end{multline}
where the constant $C$ does not depend on $x_0$, and $\alpha>0$ is defined by \eqref{eq:hyp1}. Using the classical estimate
we recalled in \eqref{eq:estimation_G}, we deduce that there exists some constant $C>0$ such that
\begin{equation}
  \label{eq:estimee_mauvaise_grad_G}
  \forall\, |x-y|>1, \quad \left|\nabla_x G(x,y)\right| \leq \frac{C}{|x-y|^{d-2}}.
\end{equation}
This estimate, applied to $H(x,y)=G(y,x),$ the Green function of the operator
$-\div\left(a^T\nabla\cdot\right)$, yields
\begin{equation}
  \label{eq:estimee_mauvaise_grad_G_2}
  \forall\, |x-y|>1, \quad \left|\nabla_y G(x,y)\right| \leq \frac{C}{|x-y|^{d-2}}.
\end{equation}
Next, we apply the proof of (i) to $\partial_{y_k}G$. More precisely, since $-\div_x(a\nabla_x\partial_{y_k}G(x,y)) = 0$ in the set
$\{|x-y|>1/2\}$, we may apply Caccioppoli inequality and the whole sequence of arguments that successively lead to
\eqref{eq:caccioppoli_appliquee} through \eqref{eq:holder0} to $\partial_{y_k} G$ instead of $G$ and we obtain, for
$1\leq q \leq 2$,
\begin{displaymath}
\int_{B_{2R}(y)\backslash B_R(y)} \left|\nabla_x \nabla_yG(x,y)\right|^q\,dx\leq {\frac
  C {R^{d(q-1)-q}}},
\end{displaymath}
where $C$ depends on $q$ but not on $y$ nor on $R$. 
In particular, summing up all these inequalities for $R = 2^k$, $k\geq 0$, we have
\begin{equation}\label{eq:32}
  \int_{|x-y|>1} |\nabla_x\nabla_yG(x,y)|^q dx \leq C\sum_{k\geq 0} \frac 1 {2^{k((d-1)q-d)}}.
\end{equation}
The right-hand side is a converging series if and only if $q>d/(d-1)$. Hence, 
\begin{equation}
  \label{eq:estimee_mauvaise_grad_grad_G}
  \forall\, q \in \left]\frac d {d-1},2\right],\quad \exists\,  C>0, \quad\forall\, y\in\RR^d,\quad \int_{\{|x-y|>1\}}|\nabla_x\nabla_y G(x,y)|^qdx\leq C.
\end{equation}
We are now going to prove that \eqref{eq:estimee_mauvaise_grad_grad_G} is valid for any $q\leq 2$, that is,
\begin{equation}
  \label{eq:estimee_mauvaise_grad_grad_G_2}
  \forall\, q \in \left]1,2\right],\quad \exists\,  C>0, \quad\forall\, y\in\RR^d,\quad \int_{\{|x-y|>1\}}|\nabla_x\nabla_y G(x,y)|^qdx\leq C,
\end{equation}
for a constant $C$ that, like in \eqref{eq:32} and \eqref{eq:estimee_mauvaise_grad_grad_G}, 
depends on $q$, but not on $y$. In order to prove \eqref{eq:estimee_mauvaise_grad_grad_G_2},
we write $G =G_{per}+ G_1+G_2$, with
\begin{equation}\label{eq:definition_Gper}
  -\div_x\left(a_0(x)\nabla_x G_{per}(x,y)\right) = \delta(x-y),
\end{equation}
\begin{equation}\label{eq:definition_G1}
  -\div_x\left(a_0(x)\nabla_x G_1(x,y)\right) = \div_x(\chi(x-y)\, b(x)\, \nabla_x G(x,y) ),
\end{equation}
\begin{equation}\label{eq:definition_G2}
  -\div_x\left(a_0(x)\nabla_x G_2(x,y)\right) = \div_x\left[(1-\chi(x-y))\, b(x)\, \nabla_x G(x,y) \right],
\end{equation}
where $\chi \in C^\infty(\RR^d)$ is a cut-off function:
\begin{displaymath}
  0\leq\chi\leq 1, \quad \chi_{|B_1(0)} = 1, \quad \chi_{|B_2(0)^c} = 0, \quad |\nabla \chi|\leq 2.
\end{displaymath}
We successively prove that $G_{per}$, $G_1$ and $G_2$ satisfy \eqref{eq:estimee_mauvaise_grad_grad_G_2}, for
$|x-y|>4$. From this we will infer that $G$ satisfies \eqref{eq:estimee_mauvaise_grad_grad_G_2}.

\noindent\textbf{Step 1: $G_{per}$ satisfies \eqref{eq:estimee_mauvaise_grad_grad_G_2}.} The results of \cite{AL1987,anantha} imply that $G_{per}$ satisfies the estimates
\begin{equation}
  \label{eq:estimations_grad_G_per}
  \forall\, (x,y) \in \RR^d\times \RR^d,\quad  |\nabla_xG_{per}(x,y)| +
  |\nabla_y G_{per}(x,y)|\leq \frac C {|x-y|^{d-1}},
\end{equation}
\begin{equation}
  \label{eq:estimations_grad_grad_G_per}
  \forall\, (x,y) \in \RR^d\times \RR^d, \quad |\nabla_x \nabla_y G_{per}(x,y) | \leq \frac C {|x-y|^{d}}.
\end{equation}
In particular, $G_{per}$ satisfies~\eqref{eq:estimee_mauvaise_grad_grad_G_2}. Actually, it even satisfies
\eqref{eq:grad-grad-G-new}.

\medskip

\noindent\textbf{Step 2: bound on $G_1$.} We prove that $G_1$ satisfies an estimate similar to
\eqref{eq:estimations_grad_grad_G_per} for $|x-y|>3$ (see \eqref{eq:estimee_G_1} below). For this purpose, we first
prove a bound on $\nabla_y G_1$, from which we deduce a bound on $\nabla_x\nabla_yG_1$. We write, from \eqref{eq:definition_G1},
\begin{displaymath}
  G_1(x,y) = \int_{\RR^d} \nabla_z G_{per}(x,z)^T b(z) \nabla_z G(z,y)\chi(z- y)dz. 
\end{displaymath}
Hence, differentiating this equality with respect to $y_k$, we have
\begin{multline}\label{eq:9}
  \partial_{y_k} G_1(x,y) = \underbrace{\int_{\RR^d} \nabla_z G_{per}(x,z)^T b(z) \nabla_z \partial_{y_k}G(z,y)\chi(z- y)dz}_{:=H_1(x,y)} \\ - \underbrace{\int_{\RR^d} \nabla_z G_{per}(x,z)^T b(z) \nabla_z G(z,y)\left(\partial_{y_k}\chi\right)(z- y)dz.}_{:=H_2(x,y)}
\end{multline}
The term $H_2(x,y)$ is easily estimated using \eqref{eq:estimee_mauvaise_grad_G_2} and \eqref{eq:estimations_grad_G_per},
and the fact that $\nabla \chi(z-y)$ vanishes outside $1<|z-y|<2$:
\begin{multline}\label{eq:10}
|H_2(x,y)|
  \\ \leq
  C\|b\|_{L^\infty(\RR^d)} \int_{1<|z-y|<2}\frac{1}{|x-z|^{d-1}} \frac 1 {|z-y|^{d-2}} dz  \leq \frac{C}{|x-y|^{d-1}},
\end{multline}
for any $x$ such that $|x-y|>3$. Next, we write $H_1(x,y)$ as follows
\begin{multline}\label{eq:16}
  H_1(x,y) = \int_{\RR^d} \left(\nabla_z G_{per}(x,z)^T b(z)-\nabla_y G_{per}(x,y)^T b(y) \right)
  \nabla_z\partial_{y_k}G(z,y) \chi(z-y)dz \\+\int_{\RR^d} \nabla_y G_{per}(x,y)^T b(y) \nabla_z\partial_{y_k}G(z,y) \chi(z-y)dz .
\end{multline}
In order to estimate the first term of the right-hand side of \eqref{eq:16}, we point out that $b\in C^{0,\alpha}_{\rm unif}(\RR^d)$, and that $G_{per}$ satisfies,
according to \cite[Theorem 3.5]{GW} and since $a_0\in C^{0,\alpha}_{\rm unif}(\RR^d)$,
\begin{multline}\label{eq:11}
  \forall\, z \text{ such that } |z-y|\leq 2, \\ \left|\nabla_z G_{per}(x,z)-\nabla_y G_{per}(x,y) \right|\leq C |z-y|^\alpha \left(\frac 1 {|x-z|^{d-1}} + \frac 1 {|x-y|^{d-1}} \right)
\end{multline}
Actually, \eqref{eq:11} is proved in \cite[Theorem 3.5]{GW} for a problem in a bounded domain with homogeneous
boundary conditions. But a careful examination of the proof shows that the constant does not depend on the
size of the domain, implying \eqref{eq:11}. In addition,
\cite[Theorem 3.3]{GW} implies that $|\nabla_z\partial_{y_k} G(z,y)|\leq C |z-y|^{-d}$ if $|z-y|\leq 2$. Hence,
using \eqref{eq:11} and the fact that $b\in \left(C^{0,\alpha}_{\rm unif}(\RR^d)\right)^d$, we infer 
\begin{multline}\label{eq:17}
  |H_1(x,y)|\leq C \int_{|z-y|<2} |z-y|^{\alpha}\left(\frac 1 {|x-z|^{d-1}} + \frac 1 {|x-y|^{d-1}} \right) \frac 1 {|y-z|^d}dz \\+ \left|\int_{\RR^d} \nabla_y G_{per}(x,y)^T b(y) \nabla_z\partial_{y_k}G(z,y) \chi(z-y)dz \right|
\end{multline}
The first term of the right-hand side of \eqref{eq:17} is dealt with using the fact that, if $|x-y|>3$ and $|y-z|<2$, then
$|x-y|\leq 3|x-z|$. Hence,
\begin{multline}\label{eq:21}
  \int_{|z-y|<2} |z-y|^{\alpha}\left(\frac 1 {|x-z|^{d-1}} + \frac 1 {|x-y|^{d-1}} \right) \frac 1 {|y-z|^d}dz\\ \leq
  \frac{1+3^{d-1}}{|x-y|^{d-1}} \int_{|z-y|<2} \frac{dz}{|y-z|^{d-\alpha}}\leq \frac{C}{|x-y|^{d-1}},
\end{multline}
where $C$ does not depend on $x$ nor on $y$. 
The second term of the right-hand side of \eqref{eq:17} is estimated using that $\chi$ has compact support, and integrating by parts:
\begin{multline}\label{eq:18}
  \left|\int_{\RR^d} \nabla_y G_{per}(x,y)^T b(y) \nabla_z\partial_{y_k}G(z,y) \chi(z-y)dz \right| \\=
  \left|\int_{\RR^d} \nabla_y G_{per}(x,y)^T b(y) \partial_{y_k}G(z,y) \nabla_z \chi(z-y)dz \right| \\
\leq \frac C {|x-y|^{d-1}} \int_{1<|z-y|<2} \frac 1 {|z-y|^{d-2}}dz\leq  \frac{C}{|x-y|^{d-1}},
\end{multline}
where $C$ is independent of $x$ and $y$.
Here, we have used \eqref{eq:estimations_grad_G_per} and \eqref{eq:estimee_mauvaise_grad_G}.
Inserting \eqref{eq:18} and~\eqref{eq:21} into \eqref{eq:17}, we have
\begin{equation}
  \label{eq:estimee_H_1}
  |H_1(x,y)|\leq \frac C {|x-y|^{d-1}}.
\end{equation}
Collecting \eqref{eq:10} and \eqref{eq:estimee_H_1}, and inserting them into \eqref{eq:9}, we find that 
\begin{equation}\label{eq:22}
  \forall\,  |x-y|>3, \quad |\nabla_y G_1(x,y)|\leq \frac C {|x-y|^{d-1}}.
\end{equation}
 
Recalling that $-\div_x(a_0\nabla_x \partial_{y_k}G_1) = 0$ in $\{x, \ |x-y|>2\}$, we may apply \cite[Lemma 16]{AL1987},
which implies that 
\begin{displaymath}
  \sup_{x\in B_R(x_0)} |\nabla_x \nabla_y G_1(x,y)| \leq \frac C R \sup_{x\in B_{2R}(x_0)} |\nabla_y G_1(x,y)|,
\end{displaymath}
for any $x_0\in \RR^d$ and $R>0$ such that $B_{2R}(x_0)\subset\{x,\ |x-y|>2\}$. Applying this with $R = \frac14 |x_0-y|-\frac34$,
we find that
\begin{equation}
  \label{eq:estimee_G_1}
  \forall\, |x_0-y|>4, \quad |\nabla_x \nabla_y G_1(x_0,y)| \leq \frac C {|x_0-y|^d}.
\end{equation}

\noindent\textbf{Step 3: $G_2$ satisfies \eqref{eq:estimee_mauvaise_grad_grad_G_2}.}
Differentiating \eqref{eq:definition_G2} with respect to $y_k$, we have
\begin{multline}\label{eq:15}
  -\div_x\left(a_0(x) \nabla_x \partial_{y_k} G_2(x,y)\right) =
  \div_x\left(b(x)\nabla_xG(x,y) \, \partial_{y_k}\left(1-\chi(x-y)\right) \right) \\+ \div_x\left(\left(1-\chi(x-y)\right)b(x) \nabla_x\partial_{y_k}
  G(x,y) \right)
\end{multline}
In the right-hand side of this equation, we notice that
\begin{equation}\label{eq:30}
  \left\|b\nabla_xG (\cdot,y)\, \partial_{y_k}\left(1-\chi(\cdot-y)\right)\right\|_{L^1\cap L^\infty(\RR^d)} \leq C,
\end{equation}
since the support of $\partial_{y_k}\left(1-\chi(x-y)\right)$ is included in $1<|x-y|<2$ and $b\nabla_x G$ is
bounded in this set. We fix an integer $n\geq 0$ such that
\begin{equation}\label{eq:25}
  (n-1)d \leq r < nd.
\end{equation}
Considering the rightmost term of \eqref{eq:15}, we have $b\in L^r\cap L^\infty(\RR^d)$ and
\eqref{eq:estimee_mauvaise_grad_grad_G}, hence, since the support of $1-\chi$ is included in $\{x\in\RR^d,\
|x-y|>1\}$,
\begin{multline}\label{eq:19}
 \left\|\left(1-\chi(\cdot-y)\right) b \nabla_x\partial_{y_k}  G(\cdot,y)\right\|_{L^{q_1}(\RR^d)}\leq C, \quad \frac 1 {q_1} = \min\left(1,\frac 1
   s + \frac 1 q\right), \\ s\in [r,+\infty], \quad q\in\left]\frac d {d-1},2\right],
\end{multline}
that is,
\begin{equation}
  \label{eq:24}
 \left\|\left(1-\chi(\cdot-y)\right) b \nabla_x\partial_{y_k}  G(\cdot,y)\right\|_{L^{q_1}(\RR^d)}\leq C, \quad \forall\, q_1\in \left]\max\left(1,\frac 1 {1 - \frac 1 d + \frac 1 r} \right),2\right],
\end{equation}
where $C$ does not depend on $y$.
Hence, \eqref{eq:15} reads
\begin{displaymath}
    -\div_x\left(a_0(x) \nabla_x \partial_{y_k} G_2\right) = \div_x(K), \quad \|K\|_{L^{q_1}}\leq C, \quad \forall\, q_1\in \left]\max\left(1,\frac 1 {1 - \frac
    1 d + \frac 1 r}\right),2\right].
\end{displaymath}
Since $a_0$ is periodic, we may apply \cite[Theorem A]{AL1991}. We thus have $\|\nabla_x\nabla_y
G_2(\cdot,y)\|_{L^{q}(\RR^d)} \leq C.$
This, together with \eqref{eq:estimations_grad_grad_G_per} and \eqref{eq:estimee_G_1}, implies that
\begin{displaymath}
  \left\|\nabla_x\nabla_y G(\cdot,y)\right\|_{L^{q_1}(\{x\in\RR^d,\ |x-y|>4\})} \leq C, \quad  \forall\, q_1\in \left]\max\left(1,\frac 1 {1 - \frac
    1 d + \frac 1 r} \right),2\right].
\end{displaymath}
In order to have this estimate on the set $\{|x-y|>1\}$ instead of $\{|x-y|>4\}$, we apply \cite[Theorem 8.32]{GT},
which implies that, since $G$ satisfies $-\div_x(a\nabla_x \partial_{y_k}G) = 0$ in $\{|x-y|>1/2\}$,
\begin{multline}\label{eq:GT-theorem 8.32-bis}
  \forall\, x_0\in \RR^d \text{ such that } |x_0-y|>1, \\ \|\partial_{y_k}G(\cdot,y)\|_{C^{1,\alpha}(B_{1/4}(x_0))} \leq C
  \|\partial_{y_k}G(\cdot,y)\|_{L^\infty(B_{1/2}(x_0))} .
\end{multline}
The right-hand side of \eqref{eq:GT-theorem 8.32-bis} is bounded using \eqref{eq:estimee_mauvaise_grad_G_2}, so we
have
\begin{displaymath}
  \left\|\nabla_x\nabla_y G(\cdot,y)\right\|_{L^{q_1}(\{x\in\RR^d,\ |x-y|>1\})} \leq C, \quad  \forall\, q_1\in \left]\max\left(1,\frac 1 {1 - \frac
    1 d + \frac 1 r} \right),2\right].
\end{displaymath}
Hence, we may repeat this argument $n$ times, where $n$ is defined by \eqref{eq:25}. Hence, we find that
\begin{equation}\label{eq:12}
  \left\|\nabla_x\nabla_y G(\cdot,y)\right\|_{L^{q_n}(\{x\in\RR^d,\ |x-y|>1\})}\leq C, \quad  \forall\, q_n\in \left]\underbrace{\max\left(1,\frac 1 {1 - \frac
    1 d + \frac n r} \right)}_{=1},2\right],
\end{equation}
where $C$ does not depend on $y$. 
Hence,
we have proved \eqref{eq:grad-grad-G-new}, but only for $q\in ]1,2]$.
In order to recover any $q>1$, we point out that, according to standard elliptic regularity
results (see for instance \cite[Theorem 7.2]{Giaquinta-Martinazzi}), 
\begin{displaymath}
  \int_{B_{1/4}(x_0)} \left|\nabla_x \partial_{y_k} G(\cdot,y)\right|^q \leq C \int_{B_{1/2}(x_0)}|\partial_{y_k} G(\cdot,y)|^q,
\end{displaymath}
where $x_0$ is such that $|x_0-y|>1$ and $C$ does not depend on $x_0$. According to
\eqref{eq:estimee_mauvaise_grad_G_2}, $\|\partial_{y_k} G(\cdot,y)\|_{L^q(\{|x-y|>1\})}\leq C$ for any $q>d/(d-2)$. Summing
all these inequalities for $x_0\in \delta \ZZ^d$, with $\delta>0$ sufficiently small, we thus have
\begin{displaymath}
  \left\|\nabla_x\nabla_y G (\cdot,y)\right\|_{L^q(\{|x-y|>1\})}\leq C, \quad \forall\, q>\frac d {d-2}.
\end{displaymath}
This finally proves \eqref{eq:grad-grad-G-new} in the case $d\geq 3$.

\medskip

We turn to the case $d=2$: here, it is not immediately clear that \eqref{eq:estimee_mauvaise_grad_G} is true,
because $G$ does not satisfy \eqref{eq:estimation_G}. However, \eqref{eq:grad-G} holds, and implies, for $q=2$, that that there exists $C>0$ such that
  \begin{equation}\label{eq:26}
    \forall\, R>0, \quad \int_{B_{2R}(y)\setminus B_R(y)} |\nabla_x G(x,y)|^2 dx \leq C.
  \end{equation}
We claim that $|\nabla_x G(x,y)|\leq C$, for all $x,y$ such that $|x-y|>1$. We prove this fact by contradiction: we assume
that there exist sequences $(y_n)_{n\in\NN}$ and $(x_n)_{n\in\NN}$ such that $|x_n-y_n|\geq 1$ and
\begin{equation}\label{eq:27}
  |\nabla_x G(x_n,y_n)| \mathop{\longrightarrow}_{n\to+\infty} +\infty.
\end{equation}
Let us define
\begin{displaymath}
  H_n(x) = \frac{G(x+x_n,y_n) - G(x_n,y_n)}{|\nabla_x G(x_n,y_n)|}.
\end{displaymath}
Then,
\begin{displaymath}
  -\div\left(a(x+x_n)\nabla H_n(x)\right) = 0,\quad \forall\, x\text{ such that }|x|<|x_n-y_n|.
\end{displaymath}
Moreover, we have
\begin{displaymath}
  H_n(0) = 0, \quad \left|\nabla H_n(0)\right| = 1.
\end{displaymath}
On the other hand,
\begin{multline}\label{eq:33}
  \int_{B_{\frac{|x_n-y_n|}3}(0)} |\nabla H_n|^2 = \frac{1}{|\nabla_x G(x_n,y_n)|^2} \int_{|x|<|x_n-y_n|/3}|\nabla_x
  G|^2(x+x_n,y_n)|dx \\
= \frac{1}{|\nabla_x G(x_n,y_n)|^2} \int_{|z-x_n|<|x_n-y_n|/3}|\nabla_x G|^2(z,y_n)|dz \\ \leq \frac{1}{|\nabla_x
  G(x_n,y_n)|^2} \int_{B_{2R_n}(y_n) \setminus B_{R_n}(y_n)}|\nabla G_x|^2(z,y_n)|dz ,
\end{multline}
where we have chosen $R_n = 2|x_n-y_n|/3.$ Indeed, this implies that 
\begin{displaymath}
\left\{z, \ |z-x_n|<|x_n-y_n|/3\right\}\subset B_{2R_n}(y_n) \setminus B_{R_n}(y_n).
\end{displaymath}
Applying \eqref{eq:26}, \eqref{eq:33} gives
\begin{displaymath}
  \int_{B_{\frac{|x_n-y_n|}3}(0)} |\nabla H_n|^2 \leq \frac C{|\nabla_x
  G(x_n,y_n)|^2} \mathop{\longrightarrow}_{n\to+\infty} 0
\end{displaymath}
Thus, $\nabla H_n\to 0$ in $\left(L^2(B_{1/3}(0))\right)^d$. Applying elliptic regularity results \cite[Theorem
8.32]{GT}, $H_n$ is bounded in $C^{1,\alpha}(B_{1/3}(0))$. Hence we may pass to the limit in the equality $|\nabla
H_n(0)| = 1$, reaching a contradiction.

Estimate
\eqref{eq:estimee_mauvaise_grad_G} is finally proved. In the present case, it reads
\begin{equation}\label{eq:28}
  \forall\, |x-y|>1, \quad |\nabla_x G(x,y)|\leq C,
\end{equation}
for some constant $C$ independent of $x$ and $y$. Since this is a point-wise estimate, we may apply it to $G(y,x)$,
which is the Green function of the adjoint operator $-\div\left(a^T\nabla\right)$. Hence, \eqref{eq:28} also holds
for $\nabla_y G$. Since $\nabla_y G$ satisfies $-\div(a\nabla_x \nabla_yG) = 0$ in the set $|x-y|>1$, we
apply here again elliptic regularity results, as for
instance \cite[Theorem 8.32]{GT}. We thus infer that we
have, instead of \eqref{eq:estimee_mauvaise_grad_grad_G}, 
\begin{equation}
  \label{eq:29}
\nabla_x\nabla_y G(x,y) \in L^\infty(\{|x-y|>1\}).
\end{equation}
Then, we adapt the proof of the case $d\geq 3$: we define here again $G_{per}$, $G_1$ and $G_2$ by
\eqref{eq:definition_Gper}-\eqref{eq:definition_G1}-\eqref{eq:definition_G2}. The first and second steps, which
deal with $G_{per}$ and $G_1$, are identical, and we do not reproduce them. The third step is different, since it
is in this step that we use \eqref{eq:estimee_mauvaise_grad_grad_G}. We write \eqref{eq:15}, and point out,
here again, that \eqref{eq:30} holds. We replace \eqref{eq:19} by the fact that
\begin{displaymath}
  \left\|(1-\chi(\cdot-y))b\nabla_x\partial_{y_k}G(\cdot,y)\right\|_{L^{q_1}(\RR^d)} \leq C, \quad \forall\, q_1\in[r,+\infty[,
\end{displaymath}
where $C$ does not depend on $y$.
Applying \cite[Theorem A]{AL1987} to \eqref{eq:15}, we infer that $\|\nabla_x\nabla_y
G_2(\cdot,y)\|_{L^{q_1}(\{|x-y|>1\})} \leq C$. Here again, this, together with \eqref{eq:estimations_grad_grad_G_per} and
\eqref{eq:estimee_G_1}, imply that 
\begin{equation}\label{eq:31}
  \left\|\nabla_x\nabla_y G(\cdot,y)\right\|_{L^{q_1}(\{|x-y|>1\})}\leq C,\quad  \forall\, q_1\in ]r,+\infty[,
\end{equation}
where $C$ does not depend on $y$.
We repeat the argument following \eqref{eq:29}, where we use \eqref{eq:31} instead of \eqref{eq:29}. This gives
\begin{displaymath}
  \left\|\nabla_x\nabla_y G(\cdot,y)\right\|_{L^{q_2}(\{|x-y|>1\})}\leq C, \quad \forall\, q_2 \in \left]\frac r 2 ,+\infty\right[
\end{displaymath}
provided $r/2\geq 1$. Otherwise we have $\|\nabla_x\nabla_y G(\cdot,y)\|_{L^{q_2}(\{|x-y|>1\})}\leq C,$ for all $q_2>1$. Repeating
the argument $n$ times, we thus have
\begin{displaymath}
  \left\|\nabla_x\nabla_y G(\cdot,y)\right\|_{L^{q_n}(\{|x-y|>1\})}\leq C, \quad \forall\, q_n\in\left]\max\left(1,\frac r n\right),+\infty\right[.
\end{displaymath}
For $n$ large enough, we thus have $\|\nabla_x\nabla_yG(\cdot,y)\|_{L^q(\{|x-y|>1\})}\leq C,$ for any $q>1$. We have
  proved \eqref{eq:grad-grad-G-new}, thereby concluding the proof of Lemma~\ref{lem:1}.\hfill $\diamondsuit$

\bigskip

\bigskip


Given Lemma~\ref{lem:1}, we now explain how one needs to modify the proof of Theorem 4.1 of \cite{cpde-defauts},
which we recall here:
\begin{theorem}[Theorem 4.1 of \cite{cpde-defauts}]
Assume that $a = a_0+b$ satisfies $0<\underline{\mu}\leq a_{0}(x)+b(x)$, $0<\underline\mu \leq a_0(x)$, a.e.,  for
some fixed constant $\underline{\mu}$, $a_{0}\in C^{0,\alpha}_{\rm unif}(\RR^d)$, $b\in C^{0,\alpha}_{\rm unif}(\RR^d)\cap L^r(\RR^d)$,
for some $r\in [1,+\infty[$. Assume that $a_0=a_{per}$ is periodic. Then,
problem (\ref{eq:correcteur}) has a  solution $w_p$ such that $w_p =
w_{p,0} + \tilde w_p$, where $w_{p,0}$ is the periodic corrector, that
is, the solution to (\ref{eq:correcteur-tilde}), and
\begin{itemize}
\item if $1\leq r< d$, then $\nabla \tilde w_p \in L^r,$ $\displaystyle \lim_{|x|\to+\infty} \tilde
  w_p(x) = 0$, and the solution $w_p$ is unique among those satisfying $w_p = v_{\rm
  per} + v,$ where $v_{\rm per}$ is periodic and $\nabla v \in L^r$;
\item if $2\leq r$, then $\nabla \tilde w_p \in L^r.$ In addition, the
  solution $w_p$ is unique in the class of solutions $w_p = v_{\rm
  per} + v,$ where $v_{\rm per}$ is periodic and $\nabla v \in L^r$.
\end{itemize}
\end{theorem}
We recall that in the proof of \cite[Theorem 4.1]{cpde-defauts}, the corrector $w_p = w_{p,per} + \tilde w_p$ is proved to exist, writing $\tilde w_p$ as 
\begin{displaymath}
  \tilde w_{p}(y)=\int \nabla_xG(y,x)\,\left[b\left(p+\nabla w_{p,per}(x)\right)\right]\,dx.
\end{displaymath}
A crucial ingredient of the proof is to establish that $\nabla \tilde w_p\in L^\infty$. 
The function $\tilde w_p$ is splitted as $\tilde w_p = w_1 + w_2$, where 
\begin{displaymath}
  w_1(y) =\int_{\RR^d}\nabla_x G(y,x)\left[b\left(p+\nabla
        w_{p,per}(x)\right)\right]\chi(x-y)\, dx, 
  \end{displaymath}
\begin{displaymath}
  w_2(y) = \int_{\RR^d}\nabla_x G(y,x)\left[b\left(p+\nabla w_{p,per}(x)\right)\right](1-\chi(x-y))\,dx,
\end{displaymath}
where, $\chi$ is a cut-off function:
$$\chi\in {\cal D}(\RR^d), \quad \chi_{|B_1(0)} = 1, \quad
\chi_{|B_2^c(0)} = 0, \quad \chi\geq 0, \quad |\nabla \chi|\leq 2.$$
Next, both $\nabla w_1$ and $\nabla w_2$ are shown to be bounded. The proof for $\nabla w_1$ is performed in
\cite{cpde-defauts}. As for $\nabla w_2$, we write
$$\nabla w_2(y) =  \int_{\RR^d}\nabla_x \nabla_y G(y,x)\left[b\left(p+\nabla
   w_{p,per}(x)\right)\right](1-\chi(x-y))\,dx,$$
thus
\begin{multline*}
  \left|\nabla w_2(y)\right| \leq \int_{|x-y|>1}|\nabla_x \nabla_y G(y,x)|\ |b(x)|\ \left(|p|+\|\nabla
    w_{p,per}\|_{L^\infty}\right) dx \\
\leq  \left(|p|+\|\nabla
    w_{p,per}\|_{L^\infty}\right) \|\nabla_x \nabla_y G(\cdot,y)\|_{L^q(\{|x-y|>1\})}\|b\|_{L^r(\RR^d)},
\end{multline*}
where we have chosen $q = r'$, that is, $\frac 1 q + \frac 1 r = 1$, and used \eqref{eq:grad-grad-G-new} in the
right-hand side. This shows that $\nabla \tilde w_p$ is bounded and the proof of Theorem~4.1 of \cite{cpde-defauts}
then proceeds unchanged.

\medskip


\begin{remark}
  We think that \eqref{eq:estimee_grad_fausse} is not true in full generality. Indeed, applying the De Giorgi-Nash estimate (see
  for example \cite[Theorem 2]{Moser}) to $|\nabla_y G|$, which is a subsolution of $-\div_x(a \nabla_x v) = 0$ in
  $\{|x-y|>0\}$, one finds
  \begin{displaymath}
    \sup_{x\in \RR^d, \ R<|x-y|<2R} |\nabla_y G(x,y)| \leq C \left(\frac 1 {|B_{2R}\setminus B_R(y)|}\int_{B_{2R}\setminus
      B_R(y)} |\nabla_y G(x,y)|^q dx\right)^{1/q},
  \end{displaymath}
where $C$ depends on $q,a$ and $d$ only. In particular it does not depend on $R$. Hence, applying
\eqref{eq:estimee_grad_fausse}, we would have
\begin{displaymath}
  \sup_{x\in \RR^d, \ R<|x-y|<2R} |\nabla_y G(x,y)| \leq \frac C {R^{d-1}},
\end{displaymath}
which would in turn imply the \emph{pointwise} estimate
\begin{equation}\label{eq:7}
  |\nabla_y G(x,y)| \leq \frac C {|x-y|^{d-1}}.
\end{equation}
This estimate is true for a periodic coefficient (see \cite{anantha,AL1991}). It is unclear for a general
coefficient. The results of \cite{anantha} give an example in which \eqref{eq:7} \emph{and} the estimate $|\nabla_x\nabla_y G(x,y)|\leq
C |x-y|^{-d}$ cannot hold \emph{together}. This is
why we do not expect \eqref{eq:estimee_grad_fausse} to hold for a general coefficient $a$. Note however that, as
stated in Remark~\ref{rk:josien} below, it is true if the coefficient $a$ is a local perturbation of a
periodic coefficient.
\end{remark}

\begin{remark}\label{rk:josien}
  For the case of a coefficient that reads $a = a^{per} + \tilde a$, with $a^{per}$ periodic and $\tilde a \in
  L^r(\RR^d)$, it is possible to adapt the proofs of \cite{AL1987,AL1989,KLS}, thereby proving directly that
  inequality \eqref{eq:7}, thus \eqref{eq:estimee_grad_fausse}, hold. The central estimate for
this is Lemma~16 of \cite{AL1987}. We will prove in \cite{josien,josien-these} that it is valid in this special case. 
\end{remark}

%
%
%

\medskip

\bibliographystyle{plain}

\end{document}